\documentclass[11pt,reqno]{amsart}
\usepackage{mathrsfs}
\usepackage{amsmath}
\usepackage{amssymb, amsmath}
\usepackage{graphicx}

\newcommand{\newcom}{\newcommand}

\newcom{\al}{\alpha}
\newcom{\be}{\beta}
\newcom{\eps}{\epsilon}
\newcom{\ga}{\gamma}
\newcom{\Ga}{\Gamma}
\newcom{\ka}{\kappa}
\newcom{\Lam}{\Lambda}
\newcom{\lam}{\lambda}
\newcom{\la}{\lambda}
\newcom{\Om}{\Omega}
\newcom{\om}{\omega}
\newcom{\Si}{\Sigma}
\newcom{\si}{\sigma}
\newcom{\tht}{\theta}
\newcom{\dtri}{\nabla}
\newcom{\tri}{\triangle}
\newcom{\oo}{\infty}
\newcom{\vphi}{\varphi}
\newcom{\cA}{{\mathcal A}}
\newcom{\cB}{{\mathcal B}}
\newcom{\cC}{{\mathcal C}}
\newcom{\cD}{{\mathcal D}}
\newcom{\cE}{{\mathcal E}}
\newcom{\cF}{{\mathcal F}}
\newcom{\cG}{{\mathcal G}}
\newcom{\cL}{{\mathcal L}}
\newcom{\cM}{{\mathcal M}}
\newcom{\cP}{{\mathcal P}}
\newcom{\cS}{{\mathcal S}}
\newcom{\cQ}{{\mathcal Q}}
\newcom{\caly}{{\mathcal Y}}
\newcom{\calZ}{{\mathcal Z}}
\newcom{\bfz}{{\bf Z}}
\newcom{\R}{\Bbb R}
\newcom{\N}{\Bbb N}
\newcom{\Z}{\Bbb Z}
\newcom{\C}{\Bbb C}
\newcom{\E}{\Bbb E}

\def\dv{\mbox{div}}

\newcom{\hn}{{\bf H}^n}
\newcom{\hnn}{{\mathbf H}^{n'}}
\newcom{\ulzs}{u^\lam_{z,s}}
\newcom{\Hl}{{{\cal  H}_\lam}}
\newcom{\fal}{F_{\al, \lam}}
\newcom{\Dh}{\Delta_{{\mathbf H}^n}}
\newcom{\fgl}{F_{\g, \lam}}

\newcom{\f}{\frac}
\newcom{\di}{\displaystyle\int}
\newcom{\ds}{\displaystyle\sum}
\newcom{\dl}{\displaystyle\lim}
\newcom{\ov}{\overline}
\newcom{\sset}{\subset}
\newcom{\wt}{\widetilde}
\newcom{\pa}{\partial}
\newcom{\p}{\partial}
\newcom\na{\nabla}
\newcom{\co}{\cdot}
\newcom{\suml}{\sum\limits}
\newcom{\supl}{\sup\limits}
\newcom{\intl}{\int\limits}
\newcom{\infl}{\inf\limits}
\newcom{\disp}{\displaystyle}
\newcom{\non}{\nonumber}
\newcom{\no}{\noindent}
\newcom{\QED}{$\square$}
\newcom{\rd}{\textrm{d}}

\def\endproof{\hphantom{MM}\hfill\llap{$\blacksquare$}\goodbreak}

\def\eqdefa{\buildrel\hbox{\footnotesize def}\over {=\!\!=}}

\newtheorem{athm}{\bf \t}[section]
\newenvironment{thm} [1] {\def\t{#1}\begin{athm} \bf \rm} {\end{athm}}
\newcom{\bthm}{\begin{thm}}\newcom{\ethm}{\end{thm}}

\newcom{\beq}{\begin{equation}}
\newcom{\eeq}{\end{equation}}
\newcom{\ben}{\begin{eqnarray}}
\newcom{\een}{\end{eqnarray}}
\newcom{\beno}{\begin{eqnarray*}}
\newcom{\eeno}{\end{eqnarray*}}
\newcom{\bali}{\begin{aligned}}
\newcom{\eali}{\end{aligned}}

\numberwithin{equation}{section}

\begin{document}

\title[Global well-posedness for the micropolar fluid system]
{Global well-posedness for the micropolar fluid system in the
critical Besov spaces}

\author{Qionglei Chen}
\address{Institute of Applied Physics and Computational Mathematics,P.O. Box 8009, Beijing 100088, P. R. China}
\email{chen\_qionglei@iapcm.ac.cn}

\author{Changxing Miao}
\address{Institute of Applied Physics and Computational Mathematics,P.O. Box 8009, Beijing 100088, P. R. China}
\email{miao\_changxing@iapcm.ac.cn}


\date{30,July 2010}



\begin{abstract}
We prove the global well-posedness for the 3-D micropolar fluid
system in the critical Besov spaces by  making a
suitable transformation to the solutions and using the Fourier
localization method, especially combined with a new $L^p$ estimate for the Green matrix to the linear system of the transformed equation.
This result allows to construct
global solutions for a class of highly oscillating initial data of Cannone's type.
Meanwhile, we analyze the long behavior of the solutions and get some decay estimates.
\end{abstract}

\maketitle

\section{Introduction}
We consider the  incompressible micropolar fluid system in
$\mathbb{R}^+\times \mathbb{R}^3$:
\begin{equation}\label{eq:micropolar}
\left\{
\begin{aligned}
&\p_tu-(\chi+\nu)\Delta u+u\cdot\na u+\na\pi-2\chi\na\times \omega=0,\\
&\p_t\omega-\mu\Delta\omega+u\cdot\na
\omega+4\chi\omega-\kappa\na\textrm{div}\omega-2\chi
\na\times u=0, \\
&{\rm div}u=0, \\
&(u,\omega)|_{t=0}=(u_0,\omega_0).
\end{aligned}
\right.
\end{equation}
Here $u(t,x)$ and $\omega(t,x)$ denote the linear velocity and the
velocity field of rotation of the fluid respectively. The scalar
$\pi(t,x)$ denotes the pressure of the fluid. The constants $\kappa,
\chi, \nu, \mu$ are the viscosity coefficients. For simplicity, we
take $\chi=\nu=\frac12$ and $\kappa=\mu=1$.

Micropolar fluid system was firstly developed by Eringen
\cite{Eringen}. It is a type of fluids which exhibits the micro-rotational effects and micro-rotational inertia, and can be  viewed as a non-Newtonian fluid. Physically,
micropolar fluid may represent fluids that consisting of rigid,
randomly oriented (or spherical particles) suspended in a viscous
medium, where the deformation of fluid particles is ignored. It can
describe many phenomena appeared in a large number of complex
fluids such as the suspensions, animal blood, liquid crystals which
cannot be characterized appropriately by the Navier-Stokes system,
and that it is important to the scientists working with the
hydrodynamic-fluid problems and phenomena. For more background, we
refer to \cite{Luka} and references therein.

If the microstructure of the fluid is not taken into account,
that is to say the effect of  the angular velocity fields of the particle's rotation is omitted,
i.e., $\omega=0$, then Eq. \eqref{eq:micropolar} reduces to the
classical Navier-Stokes equations.

Due to its importance in mathematics and physics, there is a lot of
literature devoted to  the mathematical theory of the micropolar
fluid systems. Galdi and Rionero \cite{Gald} and Lukaszewicz
\cite{Luka} proved the existence of the weak solution. The existence and uniqueness of strong solutions to the micropolar flows and the
magneto-micropolar flows either local for large data or
global for small data are considered in \cite{Boldrini,Luka,Rojas} and references therein.
Recently, inspired by the work of
Cannone and Karch \cite{Cannon-Karch} on the compressible Navier-Stokes
equations, V.-Roa and Ferreira \cite{Roa-Ferre} proved the
well-posedness of the generalized micropolar fluids system in the
pseudo-measure space which is denoted by $PM^a$-space whose Fourier
transform verifies
\begin{align}\label{eq:PM} \sup_{\R^3}|\xi|^a|\hat{f}(\xi)|<\infty.\end{align}
On the  wellposedness for the 2D case with full viscosity and partial viscosity one may refer to \cite{Luka} and \cite{Dong} respectively;
On the blow-up criterion for the smooth solution and the regularity criterion for the weak solution one refers to
\cite{yuan,Yuan} and references therein.

For the incompressible Navier-Stokes equations
\begin{align} \left\{
\begin{array}{ll}
\p_tu-\nu\Delta u+u\cdot\na u+\na p=0, \\
\dv u=0,\\
u(x,0)=u_0,
\end{array}
\right.\label{eq:NS}
\end{align}
Fujita and Kato\cite{Fuj-Kat,Kat} proved  the local
wellposedness for large initial data and the global well-posedness
for small initial data in the homogeneous Sobolev space $\dot H^{\f
12}$ and the Lebesgue space $L^3$ repectively. These spaces are all
the critical ones, which are revelent to the scaling of the
Navier-Stokes equations: if $(u,p)$ solves  (\ref{eq:NS}), then
\begin{align}\label{scaling}
(u_\lambda(t,x),p_\lambda(t,x))\eqdefa (\lambda u(\lambda^2
t,\lambda x),\lambda^2 p(\lambda^2 t,\lambda x))
\end{align}
is also a solution of (\ref{eq:NS}). The so-called {\bf critical
space} is the one such that the associated norm is invariant under
the scaling of \eqref{scaling}. Recently, Cannone \cite{Can1}(see
also \cite{Can-book}) generalized it to Besov spaces with
negative index of regularity.  More precisely, he showed that if the
initial data satisfies
$$
\|u_0\|_{\dot B^{-1+\f 3p}_{p,\infty}}\le c,\quad p>3
$$
for some small constant $c$, then the Navier-Stokes equations
(\ref{eq:NS}) is globally well-posed. Let us emphasize that this
result allows to construct global solutions for highly oscillating
initial data which may have a large norm in $\dot H^{\f 12}$ or
$L^3$. A typical example is
\begin{align*}
u_0(x)=\sin\bigl(\f {x_3} {\varepsilon}\bigr)(-\p_2\phi(x),
\p_1\phi(x),0)
\end{align*}
where $\phi\in \cS(\R^3)$ and $\varepsilon>0$ is small enough. Concerning  the compressible Navier-Stokes
equations, we have established the similar result in  the framework of the hybrid-Besov space with the help of a new estimate for
hyperbolic/parabolic system with convection terms, please refers to \cite{CMZ}. And the same idea has utilized to
 the case of the rotating Navier-Stokes equations, please refers to
\cite{CMZrotating}.

In this paper we try to prove the same result for the micropolar
fluid equation in the more natural space as the incompressible Navier-Stokes equations \eqref{eq:NS}.

Now let us sketch the main difficulty and  the strategy to overcome it.

Applying the Leray projection to the equation \eqref{eq:micropolar}, we obtain
\begin{equation}\label{eq:micropolar1}
\left\{
\begin{aligned}
&\p_tu-\Delta u+\mathbf{P}(u\cdot\na u)-\na\times \omega=0,\\
&\p_t\omega-\Delta\omega+u\cdot\na
\omega+2\omega-\na\textrm{div}\omega-
\na\times u=0, \\
&\dv u=0,\\
&(u,\omega)|_{t=0}=(u_0,\omega_0).
\end{aligned}
\right.
\end{equation}
Obviously, the system has no scaling invariant compared with the incompressible Navier-Stokes equation.
In general there are two ways to achieve  the global existence for  small  data in the critical Besov space as $\dot B^{-1+\f 3p}_{p,\infty}$
 for general $p$.
The first one is Kato's semigroup method which was extended  in
\cite{Can1}, it  turns out that the both linear terms $\na\times
\omega$ and $\na\times u$ will  play bad roles if they are regarded
as the perturbations. The second  way is to use the energy method
together with the Fourier localization technique,  but the linear
coupling effect of the system \eqref{eq:micropolar1} is too strong
to control unless the coefficients of these two linear terms are
sufficiently small,  while it is impossible.

To go around the trouble from the terms $\na\times \omega$ and $\na\times u$, we will viewed them as
certain perturbation of the Laplacian operator in some sense. More precisely, we  will  take the idea
developed in \cite{CMZ} for the compressible Navier-Stokes equations, i.e.,
investigating the following mixed linear system of Eq.\eqref{eq:micropolar1}:
\begin{equation}\label{eq:linearizedmicropolar}
\left\{
\begin{aligned}
&\p_tu-\Delta u-\na\times \omega=0,\\
&\p_t\omega-\Delta\omega+2\omega-\na\textrm{div}\omega-
\na\times u=0,
\end{aligned}
\right.
\end{equation}
and studying the action of its Green matrix which is denoted by $G(x,t)$.
From \cite{Roa-Ferre}, we have
\begin{align}\label{eq:reprensentaion of G}
\widehat{Gf}(\xi,t)=e^{-A(\xi)t}\hat{f}(\xi),
\end{align}
where
$$A(\xi)=\bigg[\begin{array}{cc}|\xi|^2I & B(\xi)\\B(\xi) & (|\xi|^2+2)I+C(\xi)\end{array}\bigg]$$with
$$B(\xi)=i\left[\begin{array}{ccc}0&-\xi_3&\xi_2\\ \xi_3&0&-\xi_1\\ -\xi_2&\xi_1&0\end{array}\right]\quad\mbox{and}\quad
C(\xi)=\left[\begin{array}{ccc}\xi_1^2&\xi_1\xi_2&\xi_1\xi_3\\ \xi_1\xi_2&\xi_2^2&\xi_2\xi_3\\ \xi_1\xi_3&\xi_2\xi_3&\xi_3^2\end{array}\right].$$
It has been shown in \cite{Roa-Ferre} that $G(x,t)$ has some similar property
with the heat kernel, i.e.,
\begin{align}\label{eq:primaryGpointwiseestimate}
\big|\widehat{G}(\xi,t)\big|\le e^{-c|\xi|^2t},
\end{align}which means that $\|{G}(x,t)f\|_{L^2}$ is bounded.
However, it is not enough  to obtain the estimates of the solution in the Besov space as we wanted.
For this purpose,  we have to analyze the behavior of the  derivative of $\widehat{G}(\xi,t)$ to
set up the  boundedness of $G(x,t)f$ in ${L^p}$. In fact,
we have the better property  which $\|{G}(x,t)f\|_{L^p}$  has exponential decay  estimate for  $\hat{f}$ supported in a ring. But
if we directly calculate its derivatives as well as utilizing the estimate \eqref{eq:primaryGpointwiseestimate}, we only have
the rough estimate for example when $\alpha=1$,
\begin{align}\label{eq:primaryGpointwiseestimateinhomo}
\big|D_\xi\widehat{G}(\xi,t)\big|\le e^{-c|\xi|^2t}t(1+|\xi|).
\end{align}
Obviously, the above is not enough for us to deduce that for any couple $(t,\lambda)$
of positive real numbers
and ${\rm supp}\, \hat{f}\subset\lambda \cC$ such that
\begin{align}\label{eq:primaryGL^pestimate}
 \|{G}(x,t)f\|_{L^p}\le Ce^{-ct\lambda^2}\|f\|_{L^p},\quad 1\le p\le\infty,
\end{align}
except that for the high frequency case $\lambda\ge1$  and for the low frequency case only $p=2$. This fact is
same to the one of the compressible Navier-Stokes equations
and the rotating
Navier-Stokes equations \cite{CMZ,CMZrotating} for which we can get the wellposedness for highly oscillating
initial data only in the  hybrid-Besov space instead  in $\dot B^{-1+\f 3p}_{p,\infty}$ as \cite{Can1}.
Owing to the speciality of the working space---the pseudo-measure space (see \eqref{eq:PM}), only the estimate \eqref{eq:primaryGpointwiseestimate}
is required in \cite{Roa-Ferre},
and their method
seems not to  work for the derivatives estimate of $\widehat{G}(\xi,t)$.

We believe that the wellposedness of \eqref{eq:micropolar} holds for highly oscillating
data  in the more natural Besov space  $\dot B^{-1+\f 3p}_{p,\infty}$ like the incompressible Navier-Stokes equation
due to the second equation of \eqref{eq:micropolar1} presents the better property
although there is  negative impact from $\na\times u$ and $\na\times \omega$.
Our ideas is to sufficiently employ  the structure properties of the systems. In
fact, we find that if making a suitable transformation to the solutions, then  Eq.  \eqref{eq:micropolar1}
reduces to a new version.
 More precisely, the vector field velocity $u=(u_1,u_2,u_3)$ is transformed  to an anti-symmetric
matrix $u_A$ with \begin{equation*}u_A\eqdefa \left(
\begin{array}{ccc}
0& u_3 & -u_2\\
-u_3& 0& u_1  \\
u_2& -u_1 & 0
\end{array}
\right),\end{equation*}
and decompose $\omega$ into
$\omega_d=\Lambda^{-1}\dv\omega$ and
$\omega_\Omega=\Lambda^{-1}\textrm{curl}\omega$, here we denote $$\Lambda^s z\eqdefa {\cF}^{-1}(|\xi|^s\widehat{z})$$ and the matrix
$$(\textrm{curl}z)_j^i\eqdefa(\partial_jz^i-\partial_iz^j)_{1\le i,j\le3}.$$
In light  of $\dv u=0$,
the system \eqref{eq:micropolar1} can be rewritten as
\begin{equation}\label{eq:micropolartransform}
\left\{
\begin{aligned}
&\p_tu_A-\Delta u_A-\Lambda \omega_\Omega=-\big(\mathbf{P}(u\cdot\na u)\big)_A,\\
&\p_t\omega_\Omega-\Delta\omega_\Omega+2\omega_\Omega-
\Lambda u_A=-\Lambda^{-1}\textrm{curl}\big(u\cdot\na\omega\big), \\
&\p_t\omega_d-2\Delta\omega_d+2\omega_d=-\Lambda^{-1}\dv\big(u\cdot\na\omega\big),\\
&\omega=\Lambda^{-1}\nabla\omega_d-\Lambda^{-1}\dv\omega_\Omega,\,\,\dv u=0, \\
&(u_A,\omega_\Omega,\omega_d)|_{t=0}=(u_{0,A},\omega_{0,\Omega},\omega_{0,d}).
\end{aligned}
\right.
\end{equation}
where $\big(\mathbf{P}(u\cdot\na u)\big)_A$ is as follows:
\begin{equation*} \left(
\begin{array}{ccc}
0& u_i\pa_i u_3-\frac{\pa_3\pa_j}{\Delta}\big(u_i\pa_i u_j\big) & -u_i\pa_i u_2+\frac{\pa_2\pa_j}{\Delta}\big(u_i\pa_i u_j\big)\vspace{.3cm}\\
-u_i\pa_i u_3+\frac{\pa_3\pa_j}{\Delta}\big(u_i\pa_i u_j\big)& 0& u_i\pa_i u_1-\frac{\pa_1\pa_j}{\Delta}\big(u_i\pa_i u_j\big) \vspace{.3cm} \\
u_i\pa_i u_2-\frac{\pa_2\pa_j}{\Delta}\big(u_i\pa_i u_j\big)& -u_i\pa_i u_1+\frac{\pa_1\pa_j}{\Delta}\big(u_i\pa_i u_j\big) & 0\vspace{.3cm}
\end{array}
\right).\end{equation*}
Let us observe the associate linear system of Eq. \eqref{eq:micropolartransform}.
Since the third equation is mainly a heat equation, we focus our attention to the first two equations of \eqref{eq:micropolartransform},
which leads us to consider the following coupling linear system:
\begin{equation}\label{eq:linearmicropolar}
\left\{
\begin{aligned}
&\p_tu_A-\Delta u_A-\Lambda\omega_\Omega=0,\\
&\p_t\omega_\Omega-\Delta\omega_\Omega+2\omega_\Omega-
\Lambda u_A=0,\\
&(u_A,\omega_\Omega)|_{t=0}=(u_{0,A},\omega_{0,\Omega}).
\end{aligned}
\right.
\end{equation}
If ${\mathcal {G}}(x,t)$ denotes  by the Green matrix of \eqref{eq:linearmicropolar}, then
$\mathcal {G}(x,t)(u_{0,A},\omega_{0,\Omega})$ is the solution of \eqref{eq:linearmicropolar}. We have
\begin{align}\label{eq:reprensentaion of cG}
\widehat{\cG f}(\xi,t)=e^{-\widetilde{A}(\xi)t}\hat{f}(\xi),
\end{align}
with $$\widetilde{A}(\xi)=\bigg[\begin{array}{cc}|\xi|^2 & |\xi|\\|\xi| & |\xi|^2+2\end{array}\bigg].$$
Then using the Laplace transform, the derivatives of  $\widehat{\cG }(\xi,t)$ can be exactly and explicitly represented, see Section  3,
which helps us to deduce
the following crucial  estimate
\begin{align*}\label{}
\big|D^{\alpha}_{\xi}\widehat{\cG}(\xi,t)\big|
\le Ce^{-c|\xi|^2t}|\xi|^{-|\alpha|}.
\end{align*}
This allows us to obtain that for any couple $(t,\lambda)$ of positive real numbers
and ${\rm supp}\, \hat{f}\subset\lambda \cC$, there holds
\begin{align*}
 \|{\cG}(x,t)f\|_{L^p}\le Ce^{-ct\lambda^2}\|f\|_{L^p},\quad 1\le p\le\infty,
\end{align*}
here $\cC$ is a ring away from zero, see Proposition \ref{Prop:Green-Lp}.
Let us emphasize that the above inequality is essential to the wellposedness in the  Besov spaces.

\bthm{Definition}
Let $1\le p\le\infty$, $T>0$. We denote $E^p_T$ by the space of functions such that
$$\|(u,\omega)\|_{E^p_{T}}\eqdefa \|(u,\omega)\|_{\widetilde{L}^\infty(0,T; \dot{B}^{\f 3p-1}_{p,\infty})}+
\|(u,\omega)\|_{\widetilde{L}^1(0,T; \dot{B}^{\f 3p+1}_{p,\infty})}<\infty.$$
If $T=\infty$, we denote $E_\infty^p$ by  $E^p$. We refer to Section 2 for the definition of $\widetilde{L}^r(X)$.

\ethm

Our main results are stated as follows.
\bthm{Theorem}\label{Thm:Besov of micropolar fluid}
There exist two positive
constants $c$ and $M$ such that for all $(u_0,\omega_0)\in \dot{B}^{\f 3p-1}_{p,\infty}$ with
\begin{align}\label{equ:initial condition}
\|u_0\|_{\dot{B}^{\f 3p-1}_{p,\infty}}+
\|\omega_0\|_{\dot{B}^{\f 3p-1}_{p,\infty}}\le c.
\end{align}
Then for $2\le p<6,$ the system (\ref{eq:micropolar}) has a global
solution $(u,\omega)\in C\big((0,\infty);\dot{B}^{\f
3p-1}_{p,\infty}\big)$ with
\begin{align*}
\|(u,\omega)\|_{L^\infty(0,\infty; \dot{B}^{\f 3p-1}_{p,\infty})}\le
M\big(\|u_0\|_{\dot B^{\f {3} p-1}_{p,\infty}}+
\|\omega_0\|_{\dot B^{\f {3} p-1}_{p,\infty}}\big).
\end{align*}
Moreover,  the uniqueness holds in $E^{p}$.
\ethm
\bthm{Remark}If we  work in the space
$\widetilde{L}^\infty(\dot{B}^{\f 3p-1}_{p,1})\cap \widetilde{L}^1(\dot{B}^{\f 3p+1}_{p,1})$, the borderline case $p=6$ can  be achieved.
Moreover, the  range of $p$ for the existence  and the uniqueness  can be extended to $[2,\infty)$ and $[2,6]$, respectively.
In fact, using the paradifferential calculus, it is easy to see that
the nonlinear term $u\cdot\na u$ and  $u\cdot\na \omega$ are bounded in $\widetilde{L}^1(\dot{B}^{\f 3p-1}_{p,1})$, i.e.,
in light of $\dv u=0$,
$$\|u\cdot\na \omega\|_{\dot{B}^{\f 3p-1}_{p,1}}\le C\|u\omega\|_{\dot{B}^{\f 3p}_{p,1}}
\le C\|u\|_{\dot{B}^{\f 3p}_{p,1}}\|\omega\|_{\dot{B}^{\f 3p}_{p,1}}, \quad\mbox{for}\quad p\in[2,\infty),$$
while  $u\omega$ is not continuous from $\dot{B}^{\f 3p}_{p,\infty}\times \dot{B}^{\f 3p}_{p,\infty}$ to $\dot{B}^{\f 3p}_{p,\infty}$.
\ethm

\bthm{Theorem}\label{Thm:C(H^12) for microploar fluid}
If $(u_0,\omega_0)\in \dot H^{\f12}$ and satisfies \eqref{equ:initial condition}, then the system
(\ref{eq:micropolar}) has a unique global solution in $C(\mathbb{R}^+; \dot H^{\f12})$.
\ethm

\bthm{Remark} Here we don't impose the $\dot{H}^{\f12}$ smallness condition on the initial data. Especially, this
allows us to obtain the global well-posedness of (\ref{eq:micropolar}) for
the highly oscillating initial velocity $(u_0,\omega_0)$. For example,
$$
u_0(x)=\sin\Big(\f {x_3} {\varepsilon}\Big)(-\p_2\phi(x),
\p_1\phi(x),0),\,\,\,\omega_0(x)=e^{i\f{x_1}{\varepsilon}}\phi(x),\quad \phi(x)\in \cS(\R^3),
$$
which satisfies
$$
\|u_0\|_{\dot B^{\f 3p-1}_{p,\infty}},\,\,\,\|\omega_0\|_{\dot B^{\f 3p-1}_{p,\infty}}\ll 1 \quad \textrm{for}\quad p>3
$$
if $\varepsilon>0$ is small enough, see Proposition \ref{Prop:oscillate}.
\ethm

Finally, we prove that the solution has the following decay estimates.
\bthm{Theorem}\label{Thm:decaytheorem}
Let $(u, \omega)$ be a solution provided by Theorem  \ref{Thm:Besov of micropolar fluid}.
Then  for all multi-indices $\alpha$, we have
\begin{align}\label{eq:decayestimate}
\|(D^\alpha_xu,D^\alpha_x\omega)\|_{\dot B^{\f3p-1}_{p,\infty}}\le C_0t^{-\f{|\alpha|}{2}},\,\,\,t>0,
\end{align}
where $C_0$ is a constant depending on the initial data.
\ethm
\bthm{Remark}From the estimate \eqref{eq:decayestimate}, one know that for $t>0$, the solution $(u, \omega)\in C^\infty(\R^3)$.

\ethm

\noindent{\bf Notation.} Throughout this paper, we  denote some notations on  the matrix $M=(M_{ij})_{1\le i,j\le m}$
$$|M|\eqdefa\sum_{i,j}|M_{ij}|,$$
and for a functional space $X$, we denote $\|M\|_X$ by
$$\|M\|_X\eqdefa\sum_{i,j}\|M_{ij}\|_X,$$

The structure of this paper is organized as follows.\vspace{.1cm}

In Section 2, we recall some basic facts about the Littlewood-Paley
theory and the functional spaces. In Section 3, we analyze Green's matrix of the
linear system \eqref{eq:linearmicropolar} and show some new
results concerning its regularizing effect.
Section 4 is devoted to the  proof of Theorem \ref{Thm:Besov of micropolar fluid}.
Section 5 is devoted to the  proof of Theorem \ref{Thm:C(H^12) for microploar fluid}.
In Section 6, we give certain decay rates of the solution.

\vspace{.2cm}

\section{Littlewood-paley theory and the function spaces}
Firstly, we introduce the Littlewood-Paley decomposition. Choose two
radial functions  $\varphi, \chi \in {\cS}(\mathbb{R}^3)$ supported in
${\cC}=\{\xi\in\mathbb{R}^3,\, \frac{3}{4}\le|\xi|\le\frac{8}{3}\}$,
${\cB}=\{\xi\in\mathbb{R}^3,\, |\xi|\le\frac{4}{3}\}$ respectively such
that
\begin{align*} \sum_{j\in\mathbb{Z}}\varphi(2^{-j}\xi)=1 \quad \textrm{for
all}\,\,\xi\neq 0.
\end{align*}
For $f\in \cS'(\R^3)$, the frequency localization operators $\Delta_j$ and $S_j(j\in\Z)$ are defined by
\begin{align}
\Delta_jf=\varphi(2^{-j}D)f,\quad S_jf=\chi(2^{-j}D)f.\nonumber
\end{align}
Moreover, we have
$$
S_jf=\sum_{k=-\infty}^{j-1}\Delta_kf\quad\textrm{ in }\quad \mathcal{Z}'(\mathbb{R}^3).
$$
Here we denote the space $\mathcal{Z}'(\mathbb{R}^3)$ by the dual
space of $\mathcal{Z}(\mathbb{R}^3)=\{f\in
{\cS}(\mathbb{R}^3);\,D^\alpha \hat{f}(0)=0;
\forall\alpha\in\big(\mathbb{N}\cup 0\big)^3
\,\mbox{multi-index}\}$.

 With our choice of $\varphi$, it is easy to verify that
\begin{align}\label{orth}
\begin{aligned}
&\Delta_j\Delta_kf=0\quad \textrm{if}\quad|j-k|\ge 2\quad
\textrm{and}
\quad \\
&\Delta_j(S_{k-1}f\Delta_k f)=0\quad \textrm{if}\quad|j-k|\ge 5.
\end{aligned}
\end{align}
For more details, please refer to \cite{Can-book,Chemin-Lecture}.

In the sequel, we will constantly use the Bony's decomposition from \cite{Bony}:
\begin{align}\label{Bonydecom}
fg=T_fg+T_gf+R(f,g), \end{align} with
$$T_fg=\sum_{j\in\mathbb{Z}}S_{j-1}f\Delta_jg,
\quad R(f,g)=\sum_{j\in\mathbb{Z}}\Delta_jf \widetilde{\Delta}_{j}g,
\quad \widetilde{\Delta}_{j}g=\sum_{|j'-j|\le1}\Delta_{j'}g.$$

Let us first recall the definition of general Besov space.

\bthm{Definition}\label{Def:Bes} Let $s\in\mathbb{R}$, $1\le p,
q\le+\infty$. The homogeneous Besov space $\dot{B}^{s}_{p,q}$ is
defined by
$$\dot{B}^{s}_{p,q}\eqdefa\big\{f\in \mathcal{Z}'(\mathbb{R}^3):\,\|f\|_{\dot{B}^{s}_{p,q}}<+\infty\big\},$$
where
\begin{align*}
\|f\|_{\dot{B}^{s}_{p,q}}\eqdefa \Bigl\|2^{ks}
\|\Delta_kf(t)\|_{L^p}\Bigr\|_{\ell^q}.\end{align*}\ethm
If  $p=q=2$,  $\dot{B}^{s}_{2,2}$ is equivalent to the homogeneous Sobelev space $\dot{H}^{s}$.
\vspace{0.1cm}

Now  let us recall Chemin-Lerner's space-time space\cite{Chemin-Lecture}.

\bthm{Definition}\label{Def:timeBes}
Let $s\in \mathbb{R},$ $1\le p, q, r\le\infty$, $I\subset
\mathbb{R}$ is an interval. The homogeneous mixed time-space Besov
space $\widetilde{L}^r(I; \dot B^s_{p,q})$ is the space of the
distribution such that
$$\widetilde{L}^r(I; \dot B^s_{p,q})\eqdefa
\{f\in {\mathcal D}(I; {\mathcal Z'}(\mathbb{R}^{d}));\,\|f\|_{\widetilde{L}^r(I; \dot{B}^{s}_{p,r})}<+\infty\}.$$
where
$$
\|f(t)\|_{\widetilde{L}^r(I;\dot B^s_{p,q})}\eqdefa
\displaystyle\bigg\|2^{sj}
\bigg(\int_{I}\|\Delta_jf(\tau)\|_{p}^r{\rm d}\tau\bigg)^{\frac{1}{r}}\bigg\|_{\ell^q(\Z)}.$$\ethm
For the convenience, we sometimes use $\widetilde{L}^r_t(\dot
B^s_{p,q})$ and $\widetilde{L}^r(\dot B^s_{p,q})$ to denote
$\widetilde{L}^r(0,t;\dot B^s_{p,q})$ and
$\widetilde{L}^r(0,\infty;\dot B^s_{p,q})$, respectively.
The direct consequence of Minkowski's inequality is that
$$L^r_t(\dot{B}^{s}_{p,q})\subseteq\widetilde{L}^r_t(\dot{B}^{s}_{p,q})
\quad\hbox{if}\quad r\le q \quad\textrm{and}\quad
\widetilde{L}^r_t(\dot{B}^{s}_{p,q})\subseteq{L}^r_t(\dot{B}^{s}_{p,q})
\quad\hbox{if}\quad r\ge q.$$

Let us state some basic properties about the Besov spaces.
\bthm{Lemma}\cite{Chemin-Lecture}\label{Lem:Besovproperty}
$(\rm{i})$ If $s<\f3p$ or $s=\f3p$ and $r=1$, then $(\dot B^s_{p,q}, \|\cdot\|_{\dot B^s_{p,q}})$ is a Banach space.\vspace{.15cm}\\
$(\rm{ii})$ We have the equivalence of norms
$$\|D^k f\|_{\dot B^s_{p,q}}\sim \|f\|_{\dot B^{s+k}_{p,q}}, \quad \textrm{for}\quad  k\in \Z^+.$$
$(\rm{iii})$ Interpolation: for
$s_1, s_2\in\R$ and $\theta\in[0,1]$, one has $$\|f\|_{\dot B^{\theta s_1+(1-\theta)s_2}_{p,q}}\le
\|f\|^\theta_{\dot B^{s_1}_{p,q}}\|f\|^{(1-\theta)}_{\dot B^{s_2}_{p,q}},$$
\ethm

The following Bernstein's lemma will be repeatedly used throughout
this paper.

\bthm{Lemma}\cite{Chemin-Lecture}\label{Lem:Bernstein}
Let $1\le p\le q\le+\infty$. Then for any $\beta,\gamma\in(\mathbb{N}\cup\{0\})^3$, there exists a constant $C$
independent of $f$, $j$ such that
\begin{align*} &{\rm supp}\hat f\subseteq
\{|\xi|\le A_02^{j}\}\Rightarrow \|\partial^\gamma f\|_{L^q}\le
C2^{j{|\gamma|}+3j(\frac{1}{p}-\frac{1}{q})}\|f\|_{L^p},
\\
&{\rm supp}\hat f\subseteq \{A_12^{j}\le|\xi|\le
A_22^{j}\}\Rightarrow \|f\|_{L^p}\le
C2^{-j|\gamma|}\sup_{|\beta|=|\gamma|}\|\partial^\beta f\|_{L^p}.
\end{align*}
\ethm
\bthm{Lemma}\cite{Cannon-Planchon}\label{Lem:positiveinequality}
Let $2\le p<+\infty$. Then for any $f$ with ${\rm supp}\hat f\subseteq \{A_12^{j}\le|\xi|\le A_22^{j}\}$.
there exists a constant $C$ independent of $f$, $j$ such that
\begin{align*}
c2^{2j}\int_{\R^3} |f|^p{\rm dx}\le
\int_{\R^3} (-\Delta f)|f|^{p-2}f{\rm dx}.
\end{align*}
\ethm
\bthm{Lemma}\cite{Chemin-Lecture}\label{Lem:paraestimate} (\rm i) Let $(s, p, r_1)$ such that $\dot B^s_
{p,r_1}$ is a Banach space. Then the paraproduct $T$
maps continuously $L^\infty\times \dot B^s_{p,r_1}$ into $\dot B^s_{p,r}$. Moreover, if $t$ is negative and $r_2$ such that
$$\f 1{r_1}+\f1{r_2}=\f1r\le 1,$$
and if $\dot B^{s+t}_{p,r}$ is a Banach space, then $T$ maps continuously $\dot B^{t}_{\infty,r_1}\times\dot B^s_{p,r_2}$ into $\dot B^{s+t}_{p,r}$.\vspace{.2cm}\\
(\rm ii) Let $(p_k, r_k)$ (for $k\in \{1, 2\}$) such that
$$s_1 +s_2>0,\quad\f1p\le\f1{p_1}+\f1{p_2}\le 1 \quad\mbox{and}\quad \f1r\le\f 1{r_1}+\f1{r_2}\le 1. $$
The operator $R$ maps $\dot B^{s_1}_{p_1,r_1}\times\dot B^{s_2}_{p_2,r_2}$
into $\dot B^{\sigma_{12}}_{p,r}$ with $$\sigma_{12}:=s_1 + s_2-3\Big(\f1{p_1}+\f1{p_2}-\f1p\Big),$$
provided that $\sigma_{12}<3/p$, or $\sigma_{12}=3/p$ and $r=1$.
\ethm

With the help of the above Lemma, we can obtain
\bthm{Lemma}\label{Lem:Product} Let $1\le p\le \infty$. Then there hold

(a)\; if $s_1, s_2\le \frac{3}{p}$ and $s_1+s_2>3\max
(0,\frac2p-1)$, then
\begin{align*}\|fg\|_{\dot B^{s_1+s_2-
\frac{3}{p}}_{p,1}} \le C\|f\|_{\dot B^{s_1}_{p,1}}\|g\|_{\dot B^{s_2}_{p,1}}.
\end{align*}
(b)\; if $s_1< \f 3 p, s_2<\f 3 p$, and $s_1+s_2> 3\max
(0,\frac2p-1)$, then
\begin{align*}\|fg\|_{\dot B^{s_1+s_2-
\frac{3}{p}}_{p,\infty}} \le C\|f\|_{\dot B^{s_1}_{p,\infty}}\|g\|_{\dot B^{s_2}_{p,\infty}}.
\end{align*}

(c)\; if $s_1\le \f 3 p, s_2<\f 3 p$, and $s_1+s_2\ge 3\max
(0,\frac2p-1)$, then
\begin{align*}\|fg\|_{\dot B^{s_1+s_2-
\frac{3}{p}}_{p,\infty}} \le C\|f\|_{\dot B^{s_1}_{p,1}}\|g\|_{\dot
B^{s_2}_{p,\infty}}.
\end{align*}
\ethm

\bthm{Proposition}\label{Prop:oscillate} Let $\phi\in \cS(\R^3)$ and $p>3$.
If $\phi_\varepsilon(x)\eqdefa e^{i\f {x_1} \varepsilon}\phi(x)$, then for any $\varepsilon>0$,
\begin{align*}
\|\phi_\varepsilon\|_{\dot B^{\f 3 p-1}_{p,\infty}}\le C\varepsilon^{1-\f 3p},
\end{align*}
here $C$ is a constant independent of $\varepsilon$.
\ethm
\noindent{\bf Proof.} Please refer to the proof of Proposition 2.9 in \cite{CMZ}, here we omit it.\endproof

\bthm{Proposition}\label{Prop:Parabolic-est}
Let $s\in\mathbb{R}$, and $p,r\in[1,\infty]$, $ \nu_1>0$, $\nu_2\ge0$.
Assume that $u_0\in \dot B^{s}_{p,q}, f\in L^1_t\dot B^{s}_{p,q}$. Then the  equation
\begin{align*}
\left\{
\begin{array}{ll}
\p_tu-\nu_1\Delta u+\nu_2 u=f,\\
u|_{t=0}=u_0,
\end{array}
\right.
\end{align*}
has a unique solution $u$ satisfying
$$\|u\|_{\widetilde{L}^r_t\dot B^{s+\f2 r}_{p,q}}
\le C\big(\|u_0\|_{\dot B^{s}_{p,q}}+
\|f\|_{ L^1_t\dot B^{s}_{p,q}}\big).$$
\ethm
\noindent{\bf Proof.} The proof is similar with the case of the heat equations, we omit it here. \endproof

\section{The linearized equations of the microfluid system}
In this section, we are devoted to analyzing the Green matrix of Eq. \eqref{eq:linearizedmicropolar}.

First, let us introduce a notation: if $(M_{ij})_{\{1\le i,j\le2\}}$ is a matrix,
$f=(f_1,f_2,f_3), g=(g_1,g_2,g_3)$ are vectors, then we denote
$$(M_{ij})_{\{1\le i,j\le2\}}\bigg(\begin{array}{ll}\!\!f \!\! \vspace{.1cm}\\ \!\!g \!\! \end{array}\bigg)\eqdefa
\bigg(\begin{array}{ll}\!\!M_{11}f+M_{12}g \!\! \vspace{.1cm}\\ \!\!M_{21}f+M_{22}g \!\! \end{array}\bigg).$$
Taking Fourier transform of \eqref{eq:linearizedmicropolar} yields that
\begin{equation}\label{eq:Fourierlinearmicropolar}
\left\{\begin{aligned}&\pa_t\widehat{u_A}+|\xi|^2\widehat{u_A}-|\xi|\widehat{\omega_\Omega}=0,\\
&\pa_t\widehat{\omega_\Omega}+(|\xi|^2+2)\widehat{\omega_\Omega}-|\xi|\widehat{u_A}=0,\\
&(\widehat{u_A},\widehat{\omega_\Omega})|_{t=0}=(\widehat{u_{0,A}},\widehat{\omega_{0,\Omega}}),
\end{aligned}\right.
\end{equation}
In what follows, we will use the   Laplace transform to get the explicit expression of $\widehat{\mathcal {G}}(\xi,t)$.

Let $p\in\sum_{\phi}$ for some $\phi\in[0,\pi/2)$, where $\sum_{\phi}=\{z\in \C \backslash \{0\},\,|\arg z|<\phi\}$.
Then we have
\begin{equation}\label{eq:LaplaceFourierlinearmicropolar}
\left\{\begin{aligned}&p(\widehat{u_A})^L+|\xi|^2(\widehat{u_A})^L-|\xi|(\widehat{\omega_\Omega})^L=\widehat{u_{0,A}},\\
&p(\widehat{\omega_\Omega})^L+(|\xi|^2+2)(\widehat{\omega_\Omega})^L-|\xi|(\widehat{u_A})^L=\widehat{\omega_{0,\Omega}}.
\end{aligned}\right.
\end{equation}
that is,
\begin{equation*}\label{}
\left(\begin{array}{ll}\!\!(\widehat{u_A})^L(\xi,t) \!\! \vspace{.3cm}\\ \!\!(\widehat{\omega_\Omega})^L(\xi,t) \!\! \end{array}\right)
=\left(\begin{array}{cc}\!\!p+|\xi|^2 & -|\xi| \vspace{.3cm}\\ \!\! -|\xi| & p+|\xi|^2+2 \!\! \end{array}\right)^{-1}
\left(\begin{array}{ll}\!\!\widehat{u_{0,A}} \!\! \vspace{.3cm}\\ \!\!\widehat{\omega_{0,\Omega}} \!\! \end{array}\right).
\end{equation*}
Setting $\lambda^2=p+|\xi|^2$, we  see that
\begin{equation*}\label{}
\left(\begin{array}{ll}\!\!(\widehat{u_A})^L \!\! \vspace{.3cm}\\ \!\!(\widehat{\omega_\Omega})^L \!\! \end{array}\right)
=\frac 1 {\det}\left(\begin{array}{cc}\!\!\lambda^2+2 & |\xi| \vspace{.3cm}\\ \!\! |\xi| & \lambda^2+2 \!\! \end{array}\right)
\left(\begin{array}{ll}\!\!\widehat{u_{0,A}} \!\! \vspace{.3cm}\\ \!\!\widehat{\omega_{0,\Omega}} \!\! \end{array}\right)
\end{equation*}
with $$\det\eqdefa\lambda^4+2\lambda^2-|\xi|^2.$$
Then we have the explicit expression of the solution of \eqref{eq:Fourierlinearmicropolar}:
\begin{equation}\label{}
\left(\begin{array}{ll}\!\!\widehat{u_A} \!\! \\ \!\!\widehat{\omega_\Omega} \!\! \end{array}\right)
=\bigg\{ \mathcal{L}^{-1}\Big(\frac{\lambda^2}{\det}\Big)I+\mathcal{L}^{-1}\Big(\frac{1}{\det}\Big)
\left(\begin{array}{cc}\!\!2&|\xi|\!\! \\ \!\!|\xi| & 0 \!\!
\end{array}\right)\bigg\}
\left(\begin{array}{ll}\!\!\widehat{u_{0,A}} \!\! \\ \!\!\widehat{\omega_{0,\Omega}} \!\! \end{array}\right),
\end{equation}
where $\mathcal{L}^{-1}$ is the reverse   Laplace transformation
with respect to $p$, and $I$ is the identity matrix. Denote
\begin{align*}
&\cA(\xi,t)\eqdefa\frac{e^{(-1-\sqrt{1+|\xi|^2})t}-e^{(-1+\sqrt{1+|\xi|^2})t}}{2\sqrt{1+|\xi|^2}},\\
&\cB(\xi,t)\eqdefa\frac{e^{(-1-\sqrt{1+|\xi|^2})t}+e^{(-1+\sqrt{1+|\xi|^2})t}}{2}.
\end{align*}
Note that
\begin{align*}
&\int_0^\infty\! e^{-pt}\big[\cA(\xi,t)+\cB(\xi,t)
\big]e^{-|\xi|^2t}{\rm d}t\\
&=\frac{\lambda^2}{(\lambda^2+1-\sqrt{1+|\xi|^2})(\lambda^2+1+\sqrt{1+|\xi|^2})}=\f{\lambda^2}{\det},
\end{align*}
and
\begin{align*}
&-\int_0^\infty\! e^{-pt}\cA(\xi,t)e^{-|\xi|^2t}{\rm d}t=
\frac{1}{(\lambda^2+1-\sqrt{1+|\xi|^2})(\lambda^2+1+\sqrt{1+|\xi|^2})}=\f{1}{\det},
\end{align*}
we obtain the following proposition.
\bthm{Proposition}\label{prop:FourierGexpression} There exists a unique solution $(\widehat{u_A},\widehat{\omega_\Omega})$ of Eq.
\eqref{eq:LaplaceFourierlinearmicropolar} which is given by
\begin{align}\label{}
\left(\begin{array}{ll}\!\!\widehat{u_A} \!\! \\ \!\!\widehat{\omega_\Omega} \!\! \end{array}\right)
&=e^{-|\xi|^2t}\Big(\widehat{\cG}_1(\xi,t)
+\widehat{\cG}_2(\xi,t)\Big)
\left(\begin{array}{ll}\!\!\widehat{u_{0,A}} \!\! \\ \!\!\widehat{\omega_{0,\Omega}} \!\! \end{array}\right)
\end{align}
with $$\widehat{\cG}_1(\xi,t)\eqdefa \cA(\xi,t) \mathcal{R}(\xi),\quad \widehat{\cG}_2(\xi,t)\eqdefa  \cB(\xi,t)I, $$where
$$\mathcal{R}(\xi)=\left(\begin{array}{cc}\!\!-1&-|\xi|\!\! \\ \!\!-|\xi| & 1 \!\!
\end{array}\right). $$
\ethm

Next we will derive the pointwise estimates for $\widehat{\cG}_1(\xi,t)$, $\widehat{\cG}_2(\xi,t)$  and their derivatives.

\bthm{Lemma}\label{lem:G1G2pointwiseestimate}For multi-indices $\alpha$, there exists a positive constant $C$ independent of $\xi$, $t$ such that
\begin{align}\label{eq:G1G2pointwiseestimate}
&|\xi|^{|\alpha|}\big|D^{\alpha}_{\xi}\widehat{\cG}_1(\xi,t)\big|,\;\;\; |\xi|^{|\alpha|}\big|D^{\alpha}_{\xi}\widehat{\cG}_2(\xi,t)\big|\nonumber\\
&\le C\Big(1+e^{\frac{|\xi|^2}{2}t}\Big)\Big((|\xi|^2 t)^{|\alpha|}+(|\xi|^2 t)^{|\alpha|-1}+\cdots+|\xi|^{2}t
+1\Big).
\end{align}
\ethm

\noindent{\bf Proof.} Mean value theorem tells us that there exists a constant $\theta\in[0,1]$ such that
$$\sqrt{1+|\xi|^2}-1=\frac{1}{2}|\xi|^2(1+|\xi|^2\theta)^{-\frac12},$$
which implies that
\begin{align}\label{eq:eestimate}
e^{(-1+\sqrt{1+|\xi|^2})t}\le e^{\frac{|\xi|^2}{2}t}.
\end{align}
Using the Leibnitz's formula yields that
\begin{align}\label{eq:Leibnitz for G_1}
D^{\alpha}_{\xi}\widehat{\cG}_1(\xi,t)=\sum_{|\alpha|=N,|\alpha_1|+|\alpha_2|=|\alpha|}&D^{\alpha_1}_{\xi}
\Big(e^{(-1-\sqrt{1+|\xi|^2})t}+e^{(-1+\sqrt{1+|\xi|^2})t}\Big)\nonumber\\&\times D^{\alpha_2}_{\xi}\left(\frac{1}{2\sqrt{1+|\xi|^2}}
\left(\begin{array}{cc}\!\!-1&-|\xi|\!\! \\ \!\!-|\xi| & 1 \!\!
\end{array}\right)\right)
\end{align}For simplicity, we only show the case of  $|\alpha|=1$ in details, the other cases ($|\alpha|>1$) can be done in the same argument. Noting that
$$1+|\xi|\le 2\sqrt{1+|\xi|^2},$$one gets
\begin{align*}
D_{\xi}\big(e^{(-1-\sqrt{1+|\xi|^2})t}\big)\frac{(1+|\xi|)}{\sqrt{1+|\xi|^2}}\le Ce^{-|\xi|^2t}t|\xi|.
\end{align*}In addition, due  to \eqref{eq:eestimate}, we obtain
\begin{align*}
D_{\xi}\big(e^{(-1+\sqrt{1+|\xi|^2})t}\big)\frac{(1+|\xi|)}{\sqrt{1+|\xi|^2}}\le Ce^{\f{|\xi|^2t}{2}}t|\xi|,
\end{align*}and
\begin{align*}
\big(&e^{(-1-\sqrt{1+|\xi|^2})t}+e^{(-1+\sqrt{1+|\xi|^2})t}\big)D_{\xi}\Big(\frac{1}{\sqrt{1+|\xi|^2}}\Big)(1+|\xi|)\nonumber\\
&\le C(e^{-|\xi|^2t}+e^{\f{|\xi|^2t}{2}})|\xi|^{-1},
\end{align*}and\begin{align*}
\big(&e^{(-1-\sqrt{1+|\xi|^2})t}+e^{(-1+\sqrt{1+|\xi|^2})t}\big)\frac{1}{\sqrt{1+|\xi|^2}}D_{\xi}
\left(\begin{array}{cc}\!\!-1&-|\xi|\!\! \\ \!\!-|\xi| & 1 \!\!
\end{array}\right)\nonumber\\& \le C(e^{-|\xi|^2t}+e^{\f{|\xi|^2t}{2}})(\sqrt{1+|\xi|^2})^{-1}\le C(e^{-|\xi|^2t}+e^{\f{|\xi|^2t}{2}})|\xi|^{-1}.
\end{align*}Combining the four above  inequalities with \eqref{eq:Leibnitz for G_1}($|\alpha|=1$), we have
\begin{align*}
|&D_{\xi}\widehat{\cG}_1(\xi,t)|\le C(1+e^{\f{|\xi|^2t}{2}})(t|\xi|+|\xi|^{-1}).
\end{align*}
Similarly, we can deduce that
\begin{align*}
|D_{\xi}\widehat{\cG}_1(\xi,t)|\le C\Big(1+e^{\frac{|\xi|^2}{2}t}\Big)\Big(&|\xi t|^{|\alpha|}+|\xi|^{|\alpha|-2} t^{|\alpha|-1}\nonumber\\&+
|\xi|^{|\alpha|-4} t^{|\alpha|-2}+\cdots+|\xi|^{-|\alpha|+2}t
+|\xi|^{-|\alpha|}\Big),
\end{align*}
from which  the estimate \eqref{eq:G1G2pointwiseestimate} holds. \endproof
\vspace{.2cm}

Thanks to Proposition \ref{prop:FourierGexpression},
we have
\bthm{Proposition}\label{eq:representation of G}
The  Fourier transform of
the Green matrix of Eq.\eqref{eq:linearmicropolar}--$\widehat{\mathcal {G}}(\xi,t)$ is shown to be
$$\widehat{\mathcal {\cG}}(\xi,t)=e^{-|\xi|^2t}\big(\widehat{\cG}_1(\xi,t)
+\widehat{\cG}_2(\xi,t)\big).$$
\ethm
\bthm{Lemma}\label{lem:Gpointwiseestimate}For any multi-indices $\alpha$, there exists a positive constant $C$ independent of $\xi$, $t$ such that
\begin{align}\label{eq:Gpointwiseestimate}
\big|D^{\alpha}_{\xi}\widehat{\cG}(\xi,t)\big|
\le Ce^{-\frac{1}{3}|\xi|^2t}|\xi|^{-|\alpha|}.
\end{align}
\ethm \noindent{\bf Proof.} Noting that for $c>\tilde{c}>0$, $k>0$,
we have $$e^{-c|\xi|^2t}(t|\xi|^2)^k\le e^{-\tilde{c}|\xi|^2t}.$$
Then  using the Leibniz formula, the estimate
$$\big|\pa^{\gamma}_{\xi}(e^{-|\xi|^2t})\big|\le C|\xi|^{-|\gamma|}e^{-\f  {11}{12}|\xi|^2t},$$
and Lemma \ref{lem:G1G2pointwiseestimate},  the estimate \eqref{eq:Gpointwiseestimate} follows easily by
 the explicit expression  of $\widehat{\cG}(\xi,t)$. \endproof
\vspace{.15cm}

Using this lemma, we can obtain the following smoothing effect on Green's matrix $\cG$, which will play an important role in this paper.

\bthm{Proposition}\label{Prop:Green-Lp} Let ${\cC}$ be a ring
centered  at 0 in $\R^3$. There exist  two positive constants $c$
and $C$ such that, for any real $p\in[1,\infty]$, any couple $(t,
\lambda)$ of positive real numbers such that if ${\rm
supp}\,\hat{u}\subset\lambda{\cC}$, then we have
\begin{align}\label{eq:Green-Lp}
&\|\mathcal{G}(x,t) u\|_{L^p}\le Ce^{-c\lambda^2
t}\|u\|_{L^p}.
\end{align}
\ethm
\noindent{\bf Proof.}\, We will adopt the spirit of the proof for heat operators as in \cite{Chemin-Lecture}. For the completeness, here we will present
a proof.

Let $\phi\in {\cD}(\mathbb{R}^3\setminus\{0\})$, which equals to 1 near the  ring ${\cC}$.
Set
$$g(t,x)\eqdefa (2\pi)^{-3}\int_{\mathbb{R}^3}e^{ix\cdot\xi}
\phi(\lambda^{-1}\xi)\widehat{\mathcal{G}}(\xi,t){\rm d}\xi.
$$
To prove \eqref{eq:Green-Lp}, it suffices to show
\begin{align}\label{eq:gL1}
\|g(x,t)\|_{L^1}\le Ce^{-c\lambda^2 t}.
\end{align}
Thanks to \eqref{eq:Gpointwiseestimate} and the support property of $\phi$, we infer that
\begin{align}\label{eq:gL2}
\int_{|x|\le\lam^{-1}}|g(x,t)|{\rm d}x&\le C
\int_{|x|\le\lam^{-1}}\int_{\mathbb{R}^3}|\phi(\lambda^{-1}\xi)||\widehat{\mathcal{G}}(\xi,t)|{\rm d}\xi {\rm d}x
\le Ce^{-c\lambda^2 t}.
\end{align}
Set $L_x\eqdefa \frac {x\cdot \na_\xi} {i|x|^2}$. Noting that $L_x(e^{ix\cdot\xi})=e^{ix\cdot\xi}$,
we get by integration by part  that
\begin{align*}
g(x,t)=& \int_{\mathbb{R}^3}L^{4}_x(e^{ix\cdot\xi})\phi(\lambda^{-1}\xi)
\widehat{\mathcal{G}}(\xi,t){\rm d}\xi\nonumber\\=&(-1)^{4}
\int_{\mathbb{R}^3}e^{ix\cdot\xi}(L^*_x)^{4}\big(\phi(\lambda^{-1}\xi)
\widehat{\mathcal{G}}(\xi,t)\big){\rm d}\xi.
\end{align*}
From the Leibniz formula and \eqref{eq:Gpointwiseestimate},
\begin{align}
&\big|(L^*_x)^{4}\big(\phi(\lambda^{-1}\xi)\widehat{\mathcal{G}}(\xi,t)\big)\big|\nonumber\\&
\le C|\lambda x|^{-4}\sum_{|\gamma|= 4,
|\beta|\le|\gamma|}\lam^{|\beta|}|(\na^{\gamma-\beta}\phi)(\lambda^{-1}\xi)|
e^{-\frac13|\xi|^2t}|\xi|^{-|\beta|}.\nonumber
\end{align}
Then we obtain, for any $\xi$ with $|\xi|\sim \lambda$,
\begin{align*}
\big|(L^*_x)^{4}\big(\phi(\lambda^{-1}\xi)\widehat{\mathcal{G}}(\xi,t)\big)\big|\le C
|\lambda x|^{-4}e^{-\f 1 3|\xi|^2t},
\end{align*}
which implies that
\begin{align*}
\int_{|x|\ge\f 1\lam}|g(x,t)|{\rm d}x\le C
e^{-c\lambda^2t}\lambda^3\int_{|x|\ge\f 1\lam}|\lambda
x|^{-4}{\rm d}x\le C
e^{-c\lambda^2t}.
\end{align*}
This together with \eqref{eq:gL2} gives \eqref{eq:gL1}.\endproof
\bthm{Proposition}\label{Prop:primaryGreen-Lp}  Let ${\cC}$ be a
ring  centered  at 0 in $\R^3$, $G(x,t)$ is the Green matrix of the
system \eqref{eq:linearizedmicropolar}, defined by
\eqref{eq:reprensentaion of G}. Then there exist two positive
constants $c$ and $C$ such that for  any couple $(t, \lambda)$ of
positive real numbers satisfying: if ${\rm
supp}\,\hat{u}\subset\lambda{\cC}$, then
\begin{align}\label{eq:Green-L2}
&\|{G}(x,t) u\|_{L^2}\le Ce^{-c\lambda^2t}\|u\|_{L^2}.
\end{align}
\ethm
\noindent{\bf Proof.}\, Thanks to Plancherel theorem and \eqref{eq:primaryGpointwiseestimate}, we get
\begin{align*}
\|G(x,t) u\|_{L^2}=\|\widehat{G}(\xi,t)\hat{u}(\xi)\|_{L^2}\le C
\|e^{-c|\xi|^2 t}\hat{u}(\xi)\|_2\le Ce^{-c\lambda^2 t}\|u\|_2,
\end{align*}
where  we have used  the support property of $\hat{u}(\xi)$.
\endproof

\section{Proof of Theorem \ref{Thm:Besov of micropolar fluid}}

\subsection{A priori estimate}
In this section, we will derive a priori estimate for the linear system \eqref{eq:micropolar1}.
\bthm{Proposition}\label{Prop:Prioriestimate} Let $2\le p<6$, $T>0$.
Assume that $(u,\omega)$ is a smooth solution of the system (\ref{eq:micropolar1})
on $[0,T]$, then we have
\begin{align}\label{equ:Prioriestimate}
\|(u,\omega)\|_{E^{p}_T}\le C
\Big(\|(u_0,\omega_0)\|_{E^{p}_0}
+\|(u,\omega)\|_{E^{p}_T}^2
\Big).
\end{align}
Here
$
\|(u_0,\omega_0)\|_{E_0^{p}}\eqdefa \|u_0\|_{\dot B^{ \frac{3}{p}-1}_{p,\infty}}+
\|\omega_0\|_{\dot B^{\frac{3}{p}-1}_{p,\infty}}.
$
\ethm
\noindent{\bf Proof.}\, Let us consider the following frequency localized system:
\begin{equation}\label{eq:linsys-local}
\left\{
\begin{aligned}
&\p_t\Delta_ju_A-\Delta \Delta_ju_A-\Lambda\Delta_j\omega_\Omega=\Delta_jF,\\
&\p_t\Delta_j\omega_\Omega-\Delta\Delta_j\omega_\Omega+2\Delta_j\omega_\Omega-
\Lambda \Delta_ju_A=\Delta_jH, \\
&(\Delta_ju_A,\Delta_j\omega_\Omega)|_{t=0}=(\Delta_ju_{0,A},\Delta_j\omega_{0,\Omega}).
\end{aligned}
\right.
\end{equation}
with $$F=-\big(\mathbf{P}(u\cdot\na u)\big)_A,\quad\quad H=-\Lambda^{-1}\textrm{curl}\big(u\cdot\na\omega\big) \quad\textrm{and}
\quad \dv u=0.$$
In terms of the Green matrix $\cG$, the solution of (\ref{eq:linsys-local}) can be expressed as
\begin{align}\label{eq:linsys-local-solution}
\left(\begin{array}{ll} \!\Delta_ju_A(t) \!\!\vspace{.15cm}\\ \!\Delta_j\omega_\Omega(t) \!\end{array}\right)=
\mathcal{G}(x,t)\left(\begin{array}{ll} \!\Delta_ju_{0,A} \!\vspace{.15cm}\\ \!\Delta_j\omega_{0,\Omega} \!\end{array}\right)
+\int_0^t\mathcal{G}(x,t-\tau)\left(\begin{array}{ll}\!\Delta_jF(\tau)\!\vspace{.15cm}\\
\!\Delta_jH(\tau)\!\end{array}\right){\rm d}\tau.
\end{align}
Applying Proposition \ref{Prop:Green-Lp} to the above equation to get
\begin{align}\label{eq:linsys-local-L^p}
\|\Delta_ju_A\|_{L^p}+\|\Delta_j\omega_\Omega\|_{L^p}\le&
Ce^{-c2^{2j} t}\bigl(\|\Delta_ju_A^0\|_{L^p}+\|\Delta_j\omega_\Omega^0\|_{L^p}\bigr)\nonumber\\
&+C\int_0^te^{-c2^{2j} (t-\tau)}
\bigl(\|\Delta_jF(\tau)\|_{L^p}+\|\Delta_jH(\tau)\|_{L^p}\bigr){\rm d}\tau.
\end{align}
Taking $L^r$ norm with respect to $t$ gives
\begin{align*}
\|\Delta_ju_A\|_{L^r_TL^p}+\|\Delta_j\omega_\Omega\|_{L^r_TL^p}\le C2^{-\f {2j} r}
\bigl(&\|\Delta_ju_{0,A}\|_{L^p}+\|\Delta_j\omega_{0,\Omega}\|_{L^p}
\\
&+\|\Delta_jF\|_{L^1_TL^p}+\|\Delta_jH\|_{L^1_TL^p}\bigr).
\end{align*}
Multiplying $2^{j(\f3p-1+\f 2r)}$ on both sides,  then taking supremum over $j\in\Z$, we derive
\begin{align*}
&\|u_A\|_{\widetilde{L}^r_T\dot B^{\f3p-1+\f 2r}_{p,\infty}}+\|\omega_\Omega\|_{\widetilde{L}^r_T\dot B^{\f3p-1+\f 2r}_{p,\infty}}\\&\le C
\big(\|u_{0,A}\|_{\dot B^{\f3p-1}_{p,\infty}}+\|\omega_{0,\Omega}\|_{\dot B^{\f3p-1}_{p,\infty}}
+\|F\|_{\widetilde{L}^1_T\dot B^{\f3p-1}_{p,\infty}}+\|H\|_{\widetilde{L}^1_T\dot B^{\f3p-1}_{p,\infty}}\big).
\end{align*}
According to the boundness of Riesz
transform on the homogeneous Besov space and Lemma \ref{Lem:Product}, we have
\begin{align}\label{eq:FGestimate}
&\big\|\big(\mathbf{P}(u\cdot\na u)\big)_A\big\|_{\widetilde{L}^1_T\dot B^{\f3p-1}_{p,\infty}}
\le C\big\|u\cdot\na u\big\|_{\widetilde{L}^1_T\dot B^{\f3p-1}_{p,\infty}}
\le C\|u\|_{\widetilde{L}^{4}_T\dot B^{\f3p-\f12}_{p,\infty}}\|\na u\|_{\widetilde{L}^{\frac43}_T\dot B^{\f3p-\f12}_{p,\infty}},
\nonumber\\&
\big\|\Lambda^{-1}\textrm{curl}(u\cdot\na\omega)\big\|_{\widetilde{L}^1_T\dot B^{\f3p-1}_{p,\infty}}
\le C\big\|u\cdot\na \omega\big\|_{\widetilde{L}^1_T\dot B^{\f3p-1}_{p,\infty}}
\le C\|u\|_{\widetilde{L}^{4}_T\dot B^{\f3p-\f12}_{p,\infty}}\|\na \omega\|_{\widetilde{L}^{\frac43}_T\dot B^{\f3p-\f12}_{p,\infty}}.
\end{align}
From the Proposition \ref{Prop:oscillate}, we infer that
\begin{align}\label{eq:omegadestimate}
\|\omega_d\|_{\widetilde{L}^r_T\dot B^{\f3p-1+\f 2r}_{p,\infty}}&\le C(\|\omega_0\|_{\dot B^{\f3p-1}_{p,\infty}}
+\big\|\Lambda^{-1}\textrm{div}(u\cdot\na\omega)\big\|_{\widetilde{L}^1_T\dot B^{\f3p-1}_{p,\infty}}\nonumber\\
&\le C\Big(\|\omega_0\|_{\dot B^{\f3p-1}_{p,\infty}}+
\|u\|_{\widetilde{L}^{4}_T\dot B^{\f3p-\f12}_{p,\infty}}\|\na \omega\|_{\widetilde{L}^{\frac34}_T\dot B^{\f3p-\f12}_{p,\infty}}\Big).
\end{align}
Thanks to the interpolation
\begin{align}\label{eq:interpolation}
\Big(\widetilde{L}^{\infty}_T\dot B^{\f3p-1}_{p,\infty},\,
\widetilde{L}^1_T\dot B^{\f3p+1}_{p,\infty}\Big)_{\frac34}=\widetilde{L}^{4}_T\dot B^{\f3p-\f12}_{p,\infty},\nonumber\\
\Big(\widetilde{L}^{\infty}_T\dot B^{\f3p-1}_{p,\infty},\,
\widetilde{L}^1_T\dot B^{\f3p+1}_{p,\infty}\Big)_{\frac14}=\widetilde{L}^{\frac43}_T\dot B^{\f3p+\f12}_{p,\infty},
\end{align}
which together with \eqref{eq:FGestimate}, \eqref{eq:omegadestimate} and Lemma \ref{Lem:Besovproperty} (ii) imply
\begin{align}\label{eq:total}
&\|(u_{A},\omega_\Omega,\omega_d)\|_{\widetilde{L}^r_T\dot B^{\f3p-1+\f 2r}_{p,\infty}}\nonumber\\&\le C
(\|u\|_{\widetilde{L}^\infty_T\dot B^{\f3p-1}_{p,\infty}}+\|u\|_{\widetilde{L}^1_T\dot B^{\f3p+1}_{p,\infty}})
(\|(u,\omega)\|_{\widetilde{L}^\infty_T\dot B^{\f3p-1}_{p,\infty}}+\|(u,\omega)\|_{\widetilde{L}^1_T\dot B^{\f3p+1}_{p,\infty}}).
\end{align}
On the other hand, noting that $\omega=\Lambda^{-1}\nabla\omega_d-\Lambda^{-1}\dv\omega_\Omega$ and
\begin{align*}
\|u\|_{\widetilde{L}^r_T\dot B^{\f3p-1+\f 2r}_{p,\infty}}=\sum_{i=1}^3\|u_i\|_{\widetilde{L}^r_T\dot B^{\f3p-1+\f 2r}_{p,\infty}}
\le \|u_A\|_{\widetilde{L}^r_T\dot B^{\f3p-1+\f 2r}_{p,\infty}},
\end{align*}
taking $r=\infty$ and $r=1$ in \eqref{eq:total}, then adding  up the resulting equations, we have
\begin{align*}\label{}
\|(u,\omega)\|_{E^p_T}\le
C\big(\|(u_0,\omega_0)\|_{\dot B^{\f3p-1}_{p,\infty}}+\|(u,\omega)\|^2_{E^p_T}\big).
\end{align*}
The proof is completed. \endproof

\subsection{Approximate solutions and uniform estimates.}
The construction of approximate
solutions is based on the following local existence theorem.
\bthm{Theorem}\label{Thm:local}\cite{yuan} Let $s>3/2$. Assume that $(u_0,\omega_0)\in H^s (\R^3)$ with $\dv u_0=0$, then there
is a positive time $T (\|(u_0,\omega_0)\|_{H^s} )$ such that a unique solution $(u,\omega)\in C([0,T ); H^s )\cap
C^1((0,T ); H^s )\cap C((0,T ); H^{s+2})$ of system \eqref{eq:micropolar} exists.

Moreover, if there exists an absolute constant $M>0$ such that if
\begin{align*}\lim_{\varepsilon\rightarrow0}
\sup_{j\in \Z} \int_{T-\varepsilon}^{T}\|\Delta_j(\na\times
u)\|_{\infty}{\rm d}t =\delta<M \end{align*} then $\delta=0$, and the
solution $(u,\omega)$ can be extended past time $t=T$. \ethm
Let us consider a sequence $(\phi_n)_{n\in \N}\in \mathcal{S}$ such that  $\phi_n$
is uniformly bounded with respect to $n$ and such that $\phi_n\equiv1$ in a
neighborhood of the ball $B(0, n)$.
Then for the initial data $u_0$, $\omega_0$, we can find a approximate sequence
$u_{0,n}=\phi_n(S_n u_0),$ and $\omega_{0,n}=\phi_n(S_n \omega_0)\in H^s$ such that
\begin{align}\label{eq:convergence}\lim_{n\rightarrow\infty}\|\phi_n(S_n u_0)-u_0\|_{\dot B^{\f3p-1}_{p,\infty}}=0,\quad
\lim_{n\rightarrow \infty}\|\phi_n(S_n \omega_0)-\omega_0\|_{\dot B^{\f3p-1}_{p,\infty}}=0.\end{align}
Then Theorem \ref{Thm:local} ensures that
there exists a maximal existence time $T_n>0$ such that the system (\ref{eq:micropolar1}) with the initial data $(u_{0,n}, \omega_{0,n})$ has
a unique solution $(u^n, \omega^n)$ satisfying
\begin{align*}
(u^n,\omega^n)\in C([0,T_n ); H^s )\cap
C^1((0,T_n ); H^s )\cap C((0,T_n ); H^{s+2}).
\end{align*}
On the other hand, using the definition of the Besov space and Lemma \ref{Lem:Bernstein}, it is easy to check that
\begin{align*}
(u^n,\omega^n)\in C([0,T_n);\dot B^{\f 3p-1}_{p,\infty})\cap L^1(0,T_n;\dot B^{\f 3p+1}_{p,\infty}).
\end{align*}

From \eqref{eq:convergence} and \eqref{equ:initial condition}
we find that
\begin{align*}
\|(u_{0,n},\omega_{0,n})\|_{\dot{B}^{\frac{3}{p}-1}_{p,\infty}}\le C_0\eta,
\end{align*}
for some constant $C_0$.
Given a constant $M$ to be chosen later on,
let us define
$$T^*_n\eqdefa  \sup\Big\{t\in[0,T_n); \,\,\|(u^n,\omega^n)\|_{L^\infty_{t}\dot B^{\f 3p-1}_{p,\infty}\cap
L^1_{t}\dot B^{\f 3p+1}_{p,\infty}}\le M\eta\Big\}.$$
Firstly, we claim that
\begin{align*}
T^*_n=T_n,\qquad \forall n\in\mathbb{N}.
\end{align*}
Using the continuity argument, it suffices to show that for all $n\in\mathbb{N}$,
\begin{align}\label{equ:uniform-est}
\|(u^n,\omega^n)\|_{\widetilde{L}^\infty_{t}\dot B^{\f 3p-1}_{p,\infty}\cap
\widetilde{L}^1_{t}\dot B^{\f 3p+1}_{p,\infty}}\le \f 34M\eta.
\end{align}
In fact, applying Proposition \ref{Prop:Prioriestimate} to obtain
\begin{align}\label{equ:appro-est1}
\|(u^{n},\omega^{n})&\|_{E^{p}_{T^*_k}}\le C
\big(C_0\eta+(M\eta)^2\big).
\end{align}
If we set $M=4CC_0$, and choose $\eta$ small enough such that
\begin{align*}
8C^2C_0\eta\le 1,
\end{align*}then the inequality (\ref{equ:uniform-est}) follows from (\ref{equ:appro-est1}). In conclusion,
we construct a sequence of approximate solution $(u^n,\omega^n)$ of (\ref{eq:micropolar1}) on $[0,T_n)$ satisfying
\begin{align}\label{equ:uniform-est1}
\|(u^n,\omega^n)\|_{E^{p}_{T_n}}\le M\eta,
\end{align}
for any $n\in \N$. Next, we claim that
\begin{align*}
T_n=+\infty,\quad \forall\, n\in\mathbb{N}.
\end{align*}
According to the Theorem \ref{Thm:local}, it remains to prove
$\na\times u^n\in \widetilde{L}^1_{T_n}\dot B^0_{\infty,\infty}$.
From \eqref{equ:uniform-est1} we know that $$\|\na\times u^n\|_{\widetilde{L}^1_{T_n}\dot B^{\frac 3p}_{p,\infty}}\le
\|\na u^n\|_{\widetilde{L}^1_{T_n}\dot B^{\frac 3p}_{p,\infty}}\le M\eta, $$
this combined with the embedding $\widetilde{L}^1_{T_n}\dot B^{\frac 3p}_{p,\infty}\hookrightarrow \widetilde{L}^1_{T_n}\dot B^0_{\infty,\infty}$
implies that $\na\times u^n\in \widetilde{L}^1_{T_n}\dot B^0_{\infty,\infty}$,
thus the continuation criterion in Theorem  \ref{Thm:local}   has been verified.

\subsection{Existence}

We will use the compact argument to prove the existence of the solution. Due to (\ref{equ:uniform-est1}),
it is easy to see that
\begin{align*}
\bullet\,\, u^n, \omega^n\textrm{ is uniformly bounded in }\widetilde{L}^\infty(0,\infty; \dot B^{\f 3p-1}_{p,\infty})
\cap \widetilde{L}^1(0,\infty;\dot B^{\f 3p+1}_{p,\infty});
\end{align*}

Let $u^n_L$, $\omega^n_L$ be a solution of
\begin{equation*}\label{}
\left\{
\begin{aligned}
&\p_t u^n_L-\Delta u^n_L=0,\quad u^n_L(0)=v_{0,n},\\
&\p_t \omega^n_L-\Delta \omega^n_L+2\omega^n_L=0,\quad \omega^n_L(0)=\omega_{0,n}.
\end{aligned}\right.\end{equation*}
It is easy to verify that $u^n_L$, $\omega^n_L$ tends to the solution of
\begin{equation}\label{eq:heatmicropolar}
\left\{
\begin{aligned}
&\p_t u_L-\Delta u_L=0,\quad u_L(0)=u_{0},\\
&\p_t \omega_L-\Delta \omega_L+2\omega_L=0,\quad \omega_L(0)=\omega_{0}.
\end{aligned}\right.\end{equation}
in $\widetilde{L}^\infty(0,\infty; \dot B^{\f 3p-1}_{p,\infty})\cap \widetilde{L}^1(0,\infty;\dot B^{\f 3p+1}_{p,\infty})$.

We set $\widetilde{u}^n\eqdefa u^n-u_{L}^n$ and $\widetilde{\omega}^n\eqdefa \omega^n-\omega^n_{L}$.
Firstly, we claim that $(\widetilde{u}^n, \widetilde{\omega}^n)$ is uniformly bounded in $C^{\f12}_{loc}(\R^+;\dot B^{\f 3 p-2}_{p,\infty})\times
C^{\f12}_{loc}(\R^+;\dot B^{\f 3 p-1}_{p,\infty}+\dot B^{\f 3 p-2}_{p,\infty})$. In fact, let us
recall that $$\p_t\widetilde{u}^n=\Delta \widetilde{u}^n-\mathbf{P}(u^n\cdot\na u^n)-\na\times \omega^n.$$
Thanks to Lemma \ref{Lem:Product}, we have
\begin{align*}\|&\mathbf{P}(u^nu^n)\|_{L^{2}\dot B^{\f 3p-1}_{p,\infty}}
\le C\|u^n\|_{L^4\dot B^{\f 3p-\f12}_{p,\infty}}\|u^n\|_{L^{4}\dot B^{\f 3p-\f12}_{p,\infty}},
\end{align*}combined with $\Delta \widetilde{u}^n\in \widetilde{L}^2(\R^+;\dot B^{\f 3 p-2}_{p,\infty})$ and
$\na\times \omega^n\in L^\infty(\R^+;\dot B^{\f 3 p-2}_{p,\infty})$  implies
$\p_t\widetilde{u}^n\in \widetilde{L}^2_{loc}(\R^+;\dot B^{\f 3 p-2}_{p,\infty})$, thus
$\widetilde{u}^n$ is uniformly bounded in $C^{\f12}_{loc}(\R^+;\dot B^{\f 3 p-2}_{p,\infty})$. On the other hand, since
$$\p_t\widetilde{\omega}^n=\Delta \widetilde{\omega}^n-2\widetilde{\omega}^n-u^n\cdot\na \omega^n-\na\times u^n,$$
by the same argument as used in the proof of $\p_t\widetilde{u}^n$,
we get $\p_t\widetilde{\omega}^n\in \widetilde{L}^2_{loc}(\R^+;\dot
B^{\f 3 p-1}_{p,\infty}+\dot B^{\f 3 p-2}_{p,\infty})$, which
implies $\widetilde{u}^n$ is uniformly bounded in
$C^{\f12}_{loc}(\R^+;\dot B^{\f 3 p-1}_{p,\infty}+\dot B^{\f 3
p-2}_{p,\infty})$. \vspace{.25cm}

Let $\{\chi_j\}_{j\in \N}$ be a sequence of smooth functions supported in the ball $B(0,j+1)$
and equal to 1 on $B(0,j)$. The
claim ensures that for any $j\in \N$, $\{\chi_j \widetilde{u}^n\}_{n\in
\N}$ is uniformly bounded in $C^\f12_{loc}(\R^+;\dot B^{
\frac{3}{p}-2}_{p,\infty})$, and $\{\chi_j\widetilde{\omega}^n\}_{k\in \N}$ is uniformly
bounded in $C^{\f{1}{2}}_{loc}(\R^+;\dot B^{\f 3 p-1}_{p,\infty}+\dot B^{\f 3 p-2}_{p,\infty})$.
Observe that for
any $\chi\in C_0^\infty(\R^+\times\mathbb{R}^3)$, for $\varepsilon\in(0,1)$, the map:
$(\widetilde{u}^n, \widetilde{\omega}^n)\mapsto  (\chi \widetilde{u}^n, \chi \widetilde{\omega}^n)$ is
compact from $$\Big(\dot B^{\f 3 p-2}_{p,\infty}\cap \dot B^{\f 3 p-1-\varepsilon}_{p,\infty}\Big)
\times\Big(\big(\dot B^{\f 3 p-1}_{p,\infty}+\dot B^{\f 3 p-1-\varepsilon}_{p,\infty}\big)\Big)\quad\hbox{into}
\quad \dot B^{\f 3 p-2}_{p,\infty}\times \big(\dot B^{\f 3 p-1}_{p,\infty}+\dot B^{\f 3 p-2}_{p,\infty}\big),$$
see \cite{Danchin-book}.
By applying Ascoli's theorem and Cantor's diagonal process,
there exists some distribution $(\widetilde{u},\widetilde{\omega})\in L^\infty\dot B^{\f 3p-1}_{p,\infty}\cap L^1\dot B^{\f 3p+1}_{p,\infty}$
such that for
any $j\in\N$, \begin{align}\label{eq:conver1}
\begin{split}
&\chi_j \widetilde{u}^n\longrightarrow \chi_j \widetilde{u} \quad \textrm{in}\quad C_{loc}(\R^+;\dot B^{\f 3 p-2}_{p,\infty}),\\
&\chi_j \widetilde{\omega}^n\longrightarrow \chi_j \widetilde{\omega} \quad \textrm{in}\quad C_{loc}
(\R^+;\dot B^{\f 3 p-1}_{p,\infty}+\dot B^{\f 3 p-2}_{p,\infty}),
\end{split}
\end{align}
With (\ref{eq:conver1}),  it is a routine process to
verify that $(\widetilde{u}+u_L,\widetilde{\omega}+\omega_L)$ satisfies the system (\ref{eq:micropolar1}) in the
sense of distribution.

Here we show as an example the case of the term $u^n\cdot\na u^n$. Let $\psi\in C_0^\infty(\R^+\times\R^3)$ and $j\in\N$ such
that $\textrm{supp} \psi\subset[0,j]\times B(0,j)$.
we write $$u^n\cdot\na u^n-u\cdot\na u=(u^n-u)\cdot\na u^n+u\cdot\na(u^n-u).$$
We will only give the estimate of the first term with help of  Bony's decomposition, and the similar argument can be applied to the term $u\cdot\na(u^n-u)$.
Thanks to $\dv u=0$ and Lemma \ref{Lem:paraestimate},
\begin{align*}
&\|T_{u^n-u}u^n\|_{L^\infty\dot B^{\f3p-3}_{p,\infty}}+\|T_{u^n}(u^n-u)\|_{L^\infty\dot B^{\f3p-3}_{p,\infty}}\\
&\le C\|{u^n-u}\|_{L^\infty\dot B^{-2}_{\infty,\infty}}\|u^n\|_{L^\infty\dot B^{\f3p-1}_{p,\infty}}
+C\|u^n\|_{L^\infty\dot B^{-1}_{\infty,\infty}}\|u^n-u\|_{L^\infty\dot B^{\f3p-2}_{p,\infty}}\\
&\le C\|{u^n-u}\|_{L^\infty\dot B^{\f3p-2}_{p,\infty}}\|u^n\|_{L^\infty\dot B^{\f3p-1}_{p,\infty}},
\end{align*}where in the last inequality we have used the embedding $\dot B^{s_1}_{p,\infty}\subseteq \dot B^{s_2}_{\infty,\infty}$ for
$s_1-\f3p=s_2$.
And
\begin{align*}
\|R(u^n-u,u^n)\|_{L^1\dot B^{\f3p-1}_{p,\infty}}
&\le C\|{u^n-u}\|_{L^\infty\dot B^{\f3p-2}_{p,\infty}}\|u^n\|_{L^1\dot B^{\f3p+1}_{p,\infty}}.
\end{align*}
The other nonlinear terms can be treated in the same way.

\subsection{Uniqueness}

In this subsection, we prove the uniqueness of the solution. Assume that
$(u^1,\omega^1)\in E^p_{T}$ and $(u^2,\omega^2)\in E^p_{T}$ are  two solutions of the system (\ref{eq:micropolar})
with the same initial data. Then we have $(\delta u,\delta \omega)=(u^1-u^2,\omega^1-\omega^2)$ satisfies
\begin{align}\label{eq:linearmicrofluid-diffe} \left\{
\begin{array}{ll}
\p_t\delta u-\Delta\delta u=\delta F,\\
\p_t\delta \omega-\Delta\delta \omega-\nabla\dv\delta\omega+2\delta\omega=\delta H, \\
(\delta a,\delta v)|_{t=0}=(0,0),
\end{array}
\right. \end{align}
where
\begin{align*}
\delta F=&\nabla\times\delta \omega-\mathbf{P}(\delta u\cdot\na u^1)-\mathbf{P}(u^2\cdot\na\delta u),\\
\delta H=&\nabla\times \delta u-\delta u\cdot\na \omega^1-u^2\cdot \na \delta\omega.
\end{align*}
Applying Proposition \ref{Prop:Parabolic-est} to  Eq. \eqref{eq:linearmicrofluid-diffe}, one obtains
\begin{align}\label{eq:velo-diff2}
&\|(\delta u(t),\delta\omega(t))\|_{\widetilde{L}^1_t\dot B^{\frac 3p}_{p,\infty}}
+\|(\delta u(t),\delta\omega(t))\|_{\widetilde{L}^2_t\dot B^{\frac 3p-1}_{p,\infty}}
\nonumber\\&\le C\|(\delta F(\tau),\delta H(\tau))\|_{\widetilde{L}^1_t\dot B^{\frac 3p-2}_{p,\infty}}.
\end{align}
From Lemma \ref{Lem:Product} and  $\dv u=0$, we infer that
\begin{align*}
&\|\delta F\|_{\widetilde{L}^1_t\dot B^{\frac 3p-2}_{p,\infty}}
+\|\delta G\|_{\widetilde{L}^1_t\dot B^{\frac 3p-2}_{p,\infty}}
\nonumber\\&\le C\|\delta u\|_{\widetilde{L}^{\f43}_t\dot B^{\frac 3p-\f12}_{p,\infty}}
\|(\omega^1,u^1,u^2)\|_{\widetilde{L}^4_t\dot B^{\frac 3p-\f12}_{p,\infty}}
+C\|(\delta \omega, \delta u)\|_{\widetilde{L}^2_t\dot B^{\frac 3p-1}_{p,\infty}}
t^{\f12}.
\end{align*}
Then we have
\begin{align}\label{eq:velo-diff2}
&\|(\delta u(t), \delta\omega(t))\|_{\widetilde{L}^1_t\dot B^{\frac 3p}_{p,\infty}}
+\|(\delta u(t), \delta\omega(t))\|_{\widetilde{L}^2_t\dot B^{\frac 3p-1}_{p,\infty}}
\nonumber\\&\le C\Big(\|\delta u\|_{\widetilde{L}^1_t\dot B^{\frac 3p}_{p,\infty}}+\|\delta u\|_{\widetilde{L}^2_t\dot B^{\frac 3p-1}_{p,\infty}}\Big)
\|(\omega^1,u^1,u^2)\|_{\widetilde{L}^2_t\dot B^{\frac 3p}_{p,\infty}}^{\f12}\|(\omega^1,u^1,u^2)\|_{\widetilde{L}^\infty_t\dot B^{\frac 3p-1}_{p,\infty}}^{\f12}
\nonumber\\&\qquad+Ct^{\f12}\|(\delta \omega,\delta u)\|_{\widetilde{L}^2_t\dot B^{\frac 3p-1}_{p,\infty}}.
\end{align}
If  $t$ is taken small enough such that $\|(\omega^1,u^1,u^2)\|_{\widetilde{L}^2_t\dot B^{\frac 3p}_{p,\infty}}$
and $t^{\f12}$ sufficiently small, then we conclude that $(\delta u,\delta \omega)=0$ on $[0,T]$, and
a continuity argument ensures that $(u^1,\omega^1)=(u^2,\omega^2)$ on $[0,\infty)$.

\section{The proof of Theorem \ref{Thm:C(H^12) for microploar fluid}}
To prove Theorem \ref{Thm:C(H^12) for microploar fluid} , we will use the Green matrix of the linear system \eqref{eq:linearizedmicropolar}.
Let us return to \eqref{eq:micropolar1}. Due to $\dv u=0$, we have
\begin{align}\label{eq:Green-solution}
\bigg(\begin{array}{ll} \!u\!\!\vspace{.1cm}\\ \!\omega \!\end{array}\bigg)&=
{G}(x,t)\left(\begin{array}{ll} \!u_0 \!\vspace{.1cm}\\ \!\omega_0 \!\end{array}\right)
-\int_0^t{G}(x,t-\tau)\na \cdot\left(\begin{array}{ll}\!\!\mathbf{P}(uu)\!\!\vspace{.1cm}\\
\!\!u \omega\!\!\end{array}\right){\rm d}\tau\nonumber\\ &=\left(\begin{array}{ll} \!G_{ij}(t)u_0^j \!\vspace{.1cm}\\ \!G_{(i+3)j}(t)\omega_0^j \!\end{array}\right)
-\int_0^t\left(\begin{array}{ll}\!\!G_{ij}\pa_k\mathbf{P}(u_ku_j)+G_{i(j+3)}\pa_k(u_k\omega_j)\!\!\vspace{.1cm}\\
\!\!G_{(i+3)j}\pa_k\mathbf{P}(u_ku_j)+G_{(i+3)(j+3)}\pa_k(u_k\omega_j)\!\!\end{array}\right){\rm d}\tau
\nonumber\\ &\eqdefa {G}(t)(u_0,\omega_0)
+\left(\begin{array}{ll}\!B(u,\omega)\!\vspace{.15cm}\\
\!\widetilde{B}(u,\omega)\!\end{array}\right), \quad i=1,2,3,
\end{align}
here $G_{ij}(x,t)$ is the element of the Green matrix $G(x,t)$,
and the summation convention over
repeated indices $1\le j,k\le3$ is used.

In view of  the relationship:
$\dot{H}^{\f12}\thickapprox\dot B^{\f12}_{2,2}$,  we have
\begin{align*}\label{}
&\|B(u,\omega)\|_{\widetilde{L}^\infty_T\dot{H}^{\f12}}\le
\Big\|\int_0^t G(t-\tau)\na \cdot \big(\mathbf{P}(uu)+u  \omega\big)(\tau)d\tau\Big\|_{\widetilde{L}^\infty_T\dot{H}^{\f12}}\nonumber\\
&\le C\Big(\sum_{j\in\Z}2^j\Big(\sup_{t\in[0,T)}\int_0^t\|G(t-\tau)\na \cdot \Delta_j(u u+u \omega)(\tau)\|_{L^2}d\tau\Big)^2\Big)^{\f12}
\nonumber\\&\le C
\Big\|2^{\f32j}\sup_{t\in[0,T)}\int_0^te^{-c2^{2j}t}\|\Delta_j(u  u+u  \omega)(\tau)\|_{L^2}d\tau\Big\|_{\ell^2}\nonumber\\&\le
C\Big(\|T_uu\|_{\widetilde{L}^4_T\dot B^{0}_{2,2}}+\|T_u\omega+T_\omega u\|_{\widetilde{L}^4_T\dot B^{0}_{2,2}}+
\|R(u,u)+R(u,\omega)\|_{\widetilde{L}^{\f43}_T\dot B^{1}_{2,2}}\Big),
\end{align*}where in the third inequality we have used used Lemma \ref{Lem:Bernstein} and Proposition \ref{Prop:primaryGreen-Lp},
in the last inequality we have used  Bony's  decomposition.
From Lemma \ref{Lem:paraestimate}, we have
\begin{align*}
&\|T_u\omega\|_{\widetilde{L}^4_T\dot B^{0}_{2,2}}
\le C\|u\|_{\widetilde{L}^4_T\dot B^{-\f12}_{\infty,\infty}}
\|\omega\|_{\widetilde{L}^{\infty}_T\dot B^{\f12}_{2,2}}\le C\|u\|_{\widetilde{L}^4_T\dot B^{\f3p-\f12}_{p,\infty}}
\|\omega\|_{\widetilde{L}^{\infty}_T\dot B^{\f12}_{2,2}},
\\
&\|T_{\omega}u\|_{\widetilde{L}^4_T\dot B^{0}_{2,2}}
\le C\|\omega\|_{\widetilde{L}^{4}_T\dot B^{-\f12}_{\infty,\infty}}
\|u\|_{\widetilde{L}^\infty_T\dot B^{\f12}_{2,2}}\le C\|\omega\|_{\widetilde{L}^{4}_T\dot B^{\f3p-\f12}_{p,\infty}}
\|u\|_{\widetilde{L}^\infty_T\dot B^{\f12}_{2,2}},
\end{align*}
and
\begin{align*}
\|R(u,\omega)\|_{\widetilde{L}^{\f43}_T\dot B^1_{2,2}}&\le C
\|u\|_{\widetilde{L}^{\f43}_T\dot B^{\f12}_{\infty,\infty}}
\|\omega\|_{\widetilde{L}^{\infty}_T\dot B^{\f12}_{2,2}}\le C\|u\|_{\widetilde{L}^{\f43}_T\dot B^{\f3p+\f12}_{p,\infty}}
\|\omega\|_{\widetilde{L}^{\infty}_T\dot B^{\f12}_{2,2}}.
\end{align*}
The terms $T_uu$ and $R(u,u)$ can be treated in the same way as  $T_u\omega$,  $R(u,\omega)$, respectively.
Combining  the  above inequalities, we obtain
\begin{align}\label{equ:reg:binear term2}
\|B(u,\omega)\|_{\widetilde{L}^\infty_T\dot H^{\f12}}\le C\|(u,\omega)\|_{E^p_{T}}\|(u,\omega)\|_{\widetilde{L}^\infty_T\dot H^{\f12}}.
\end{align}
Similarly, we have
\begin{align}\label{equ:reg:binear term1}
\|\widetilde{B}(u,\omega)\|_{\widetilde{L}^\infty_T\dot H^{\f12}}\le C\|(u,\omega)\|_{E^p_{T}}\|(u,\omega)\|_{\widetilde{L}^\infty_T\dot H^{\f12}}.
\end{align}
From Proposition \ref{Prop:primaryGreen-Lp}, it is easy to verify
that
\begin{align}\label{eq:regular:initial data}
\|G(t)(u_0,\omega_0)\|_{L^\infty_T\dot{H}^{\f12}}\le C\|e^{-c2^{2j}t}\|_{L^\infty_T}\|(u_0,\omega_0)\|_{\dot H^{\f12}}
\le C\|(u_0,\omega_0)\|_{\dot H^{\f12}}.
\end{align}
It follows from the Theorem 1.2  that $\|(u,\omega)\|_{E^p_T}\le \eta $, then if $\eta$ is sufficiently small such that
$\eta C\le\f12$, we have for any $T>0$\begin{align*}
\|(u,\omega)\|_{L^\infty_T\dot{H}^{\f12}}\le 2C\|(u_0,\omega_0)\|_{\dot H^{\f12}}.
\end{align*}This finishes the existence of the proof of the Theorem \ref{Thm:C(H^12) for microploar fluid}.

\vspace{.2cm}

\noindent{\bf The uniqueness in $C(\dot H^{\f12})$ .}\,\,We will
adopt the spirit of \cite{Can-book}. Firstly, let us recall the
following bilinear estimate from \cite{CMZrotating}:
\bthm{Lemma}\label{Le:uniquebilinearestimate}For any  $T>0$, the
bilinear operators $B(u,v)(t)$, $\widetilde{B}(u,v)(t)$ are
bi-continuous from $L^\infty_T(\dot B^{\f12}_{2,\infty})\times
L^\infty_T(\dot H^{\f12})$ to $L^\infty_T(\dot
B^{\f12}_{2,\infty})$. Furthermore, we have
\begin{align*}
\|B(u,v)\|_{L^\infty_T\dot{B}^{\f12}_{2,\infty}}\,,\,\,\|\widetilde{B}(u,v)\|_{L^\infty_T\dot{B}^{\f12}_{2,\infty}}\le C
\|u\|_{\widetilde{L}^\infty_T\dot{B}^{\f12}_{2,\infty}}
\Big\|(e_{k,T})^{\f14}2^{\f k2}\|\Delta_jv\|_{L^\infty_TL^2}\Big\|_{\ell^2(k\in\Z)},
\end{align*}
here $$e_{k,T}\eqdefa 1-e^{-c2^{2k}T},$$
where $c>0$ is a constant independent of $j,T, u, v$.
\ethm

Now let $(u,\omega)$ and $(v,\varpi)$ be two solutions in $C(0,T;\dot H^{\f12})$ with the initial data $(u_0,\omega_0)\in \dot H^{\f12}$.
Using \eqref{eq:Green-solution}, we have the difference
\begin{align*}
u-v=&B\big(u-G_{ij}(t)u_{0}^j, u-v\big)+B\big(G_{ij}(t)u_{0}^j, u-v\big)
+B\big(u-v, v-G_{ij}(t)u_{0}^j\big)\nonumber\\&+B\big(u-v, G_{ij}(t)u_{0}^j\big)+B\big(u-G_{ij}(t)u_{0}^j, \omega-\varpi\big)
+B\big(G_{ij}(t)u_{0}^j, \omega-\varpi\big)\nonumber\\&+B\big(u-v, \varpi-G_{(i+3)j}(t)\omega_{0}^j\big)
+B\big(u-v, G_{(i+3)j}(t)\omega_{0}^j\big), \quad i=1,2,3.
\end{align*}
We have the same representaion for $\omega-\varpi$ replacing $B$ by $\widetilde{B}$.
We get by Lemma \ref{Le:uniquebilinearestimate} that
\begin{align*}\label{}
&\sup_{t\in(0,T)}\big(\|(u-v)(t)\|_{\dot B^{\f12}_{2,\infty}}+\|(\omega-\varpi)(t)\|_{\dot B^{\f12}_{2,\infty}}\big)\nonumber\\
&\quad\le C\sup_{t\in(0,T)}\big(\|(u-v)(t)\|_{\dot B^{\f12}_{2,\infty}}+\|(\omega-\varpi)(t)\|_{\dot B^{\f12}_{2,\infty}}\big)
\nonumber\\&\qquad\times\Big(\big\|(1-e^{-c2^{2k} T})^{\f 14}2^{\f k2}\big(\|\Delta_j u_0\|_2+\|\Delta_j \omega_0\|_2\big)\big\|_{\ell^2}
\nonumber\\&\qquad\qquad+\sup_{t\in(0,T)}(\|u-G(t)u_0\|_{\dot H^{\f12}}
+\|v-G(t)u_0\|_{\dot H^{\f12}}+\|\varpi-G(t)\omega_0\|_{\dot H^{\f12}})\Big).
\end{align*}
With the help of the fact: if $T$ is chosen sufficiently small and $(u_0,\omega_0)\in\dot H^{\f12}$, then
\begin{align*}\label{}
\big\|(1-e^{-c2^{2k} T})^{\f 14}2^{\f k2}\big(\|\Delta_k u_0\|_2+\|\Delta_k \omega_0\|_2\big)\big\|_{\ell^2}\le\f14
\end{align*}
and the strong continuity in time of the $\dot H^{\f12}$ norm of the Duhamel's term  of the solution
$(u,\omega)$ and $(v,\varpi)$,
then  a small enough time  $T$ is to be chosen such that
the three terms in the blank is dominated by 1/2 which implies that
$\|(u-v,\omega-\varpi)(t)\|_{\dot B^{\f12}_{2,\infty}}\equiv0$ on $[0,T]$. Then by the standard
argument ensures that $u=v, \omega=\varpi$ on $[0,\infty)$.

\section{The decay estimate}
Set $$W(T)\eqdefa \sup_{0\le t\le T,0<|\alpha|}t^{\frac{|\alpha|} 2}\Big(\|D^\alpha_xu\|_{\dot B^{\f3p-1}_{p,\infty}}+
\|D^\alpha_x\omega\|_{\dot B^{\f3p-1}_{p,\infty}}\Big).$$
Taking $D^\alpha_x$ on both sides of \eqref{eq:linsys-local-solution}, one gets
\begin{align*}\label{}
\left(\begin{array}{ll} \!\Delta_jD^\alpha_xu_A \!\!\vspace{.15cm}\\ \!\Delta_jD^\alpha_x\omega_\Omega \!\end{array}\right)=
D^\alpha_x\mathcal{G}(\cdot,t)\left(\begin{array}{ll} \!\Delta_ju_{0,A} \!\vspace{.15cm}\\ \!\Delta_j\omega_{0,\Omega} \!\end{array}\right)
+\int_0^tD^\alpha_x\mathcal{G}(\cdot,t-\tau)\left(\begin{array}{ll}\!\Delta_jF(\tau)\!\vspace{.15cm}\\
\!\Delta_jH(\tau)\!\end{array}\right){\rm d}\tau.
\end{align*}
Applying Lemma \ref{Lem:Bernstein} to the above equation, we have
\begin{align}\label{eq:decay}
\|\Delta_jD^\alpha_xu_A\|_{L^p}+\|\Delta_jD^\alpha_x\omega_\Omega\|_{L^p}\le &Ce^{-c2^{2j}t}2^{j|\alpha|}
\bigl(\|\Delta_ju_{0,A}\|_{L^p}+\|\Delta_j\omega_{0,\Omega}\|_{L^p}\bigr)\nonumber\\
&+\mathcal{I}+\mathcal{II}
\end{align}
where
\begin{align*}
&\mathcal{I}=C\int_0^{t/2}2^{j|\alpha|}\big(\|\mathcal{G}(t-\tau)\Delta_jF(\tau)\|_{L^p}+\|\mathcal{G}(t-\tau)\Delta_jG(\tau)\|_{L^p}\big)d\tau,\nonumber\\
&\mathcal{II}=C\int_{t/2}^t2^{j}\big(\|\mathcal{G}(t-\tau)D^{\alpha-1}_x\Delta_jF(\tau)\|_{L^p}
+\|\mathcal{G}(t-\tau)D^{\alpha-1}_x\Delta_jH(\tau)\|_{L^p}\big)d\tau.
\end{align*}
Noting that the inequality \begin{align}\label{eq:expotentialrelation}
e^{-ct2^{2j}}2^{j|\beta|}\le e^{-\tilde{c}t2^{2j}}t^{-\f {|\beta|}2},\quad |\beta|\ge0,
\end{align}
and Proposition \ref{Prop:Green-Lp}, we get that
\begin{align*}\label{}
&\mathcal{I}\le C\int_0^{t/2}e^{-\tilde{c}2^{2j} (t-\tau)}(t-\tau)^{-\frac{|\alpha|}{2}}
\bigl(\|\Delta_jF(\tau)\|_{L^p}+\|\Delta_jH(\tau)\|_{L^p}\bigr){\rm d}\tau\nonumber\\
&\mathcal{II}\le C\int_{t/2}^te^{-\tilde{c}2^{2j} (t-\tau)}(t-\tau)^{-\frac{1}{2}}
\bigl(\|D^{\alpha-1}_x\Delta_jF(\tau)\|_{L^p}+\|D^{\alpha-1}_x\Delta_jH(\tau)\|_{L^p}\bigr){\rm d}\tau
\end{align*}
In the following  we denote by $c_j(j\in\Z)$  a sequence in $\ell^1$ with the norm $\|\{c_j\}\|_{\ell^1}\le 1$.
In light of  \eqref{eq:FGestimate} and interpolation \eqref{eq:interpolation}, the straightforward calculation  shows that
\begin{align}\label{eq:decay:I}
\mathcal{I}&\le Ct^{-\frac{|\alpha|}{2}}\int_0^{t/2}e^{-\tilde{c}2^{2j} (t-\tau)}\bigl(\|\Delta_jF(\tau)\|_{L^p}+\|\Delta_jH(\tau)\|_{L^p}\bigr){\rm d}\tau
\nonumber\\&\le  Cc_j2^{-j(\f3p-1)}t^{-\frac{|\alpha|}{2}}\big(\|F\|_{\widetilde{L}^1_T\dot B^{\f3p-1}_{p,\infty}}+\|H\|_{\widetilde{L}^1_T\dot B^{\f3p-1}_{p,\infty}}\big)
\nonumber\\&\le  Cc_j2^{-j(\f3p-1)}t^{-\frac{|\alpha|}{2}}\|(u,\omega)\|^2_{E_T^p}
\nonumber\\&\le  Cc_j2^{-j(\f3p-1)}t^{-\frac{|\alpha|}{2}}\|(u_0,\omega_0)\|_{E_0^p}^2.
\end{align}
Thanks to the H\"{o}lder inequality, we have
\begin{align*}
\mathcal{II}&\le C\|e^{-c2^{2j} t}\|_{L^4_T}\Big(\int_{t/2}^t(t-\tau)^{-\frac{2}{3}}\textrm{d}\tau\Big)^{\frac{3}{4}}
\bigl(\|D^{\alpha-1}_x\Delta_jF\|_{L^\infty_TL^p}+\|D^{\alpha-1}_x\Delta_jH\|_{L^\infty_TL^p}\bigr)\nonumber\\&
\le C2^{-\f j2}t^{\f14}\bigl(\|D^{\alpha-1}_x\Delta_jF\|_{L^\infty_TL^p}+\|D^{\alpha-1}_x\Delta_jH\|_{L^\infty_TL^p}\bigr).
\end{align*}
The divergence free condition on $u$, Lemma \ref{Lem:Bernstein}  and
Lemma \ref{Lem:Product} give that
\begin{align*}
&\|\Delta_jD^{\alpha-1}_xH\|_{L^\infty_TL^p}\le C2^j\big\|\Delta_j\big((D^{\alpha-1}_x u)\omega+u(D^{\alpha-1}_x \omega)\big)\big\|_{L^\infty_TL^p}
\\&\le Cc_j2^{-j(\f3p-\f32)}\big\|(D^{\alpha-1}_x u)\omega+u(D^{\alpha-1}_x \omega)\big\|_{L^\infty_T\dot B^{\f3p-\f12}_{p,\infty}}
\\&
\le Cc_j2^{-j(\f3p-\f32)}\Big(\|D^{\alpha-1}_x u\|_{L^\infty_T\dot B^{\f3p-\f14}_{p,\infty}}
\|\omega\|_{L^\infty_T\dot B^{\f3p-\f14}_{p,\infty}}+\|D^{\alpha-1}_x \omega\|_{L^\infty_T\dot B^{\f3p-\f14}_{p,\infty}}
\|u\|_{L^\infty_T\dot B^{\f3p-\f14}_{p,\infty}}\Big).
\end{align*}
By means of Interpolation and Lemma \ref{Lem:Besovproperty} (ii), we
have
\begin{align*}
\|D^{\alpha-1}_x u\|_{L^\infty_T\dot B^{\f3p-\f14}_{p,\infty}}&\le
C\|u\|_{L^\infty_T{\dot B^{\f3p-1+\alpha}_{p,\infty}}}^{1-\f1{4\alpha}}\|u\|_{L^\infty_T{\dot B^{\f3p-1}_{p,\infty}}}^{\f1{4\alpha}}
\le C\|D^{\alpha}_x u\|_{L^\infty_T{\dot B^{\f3p-1}_{p,\infty}}}^{1-\f1{4\alpha}}\|u\|_{L^\infty_T{\dot B^{\f3p-1}_{p,\infty}}}^{\f1{4\alpha}}
\\&\le Ct^{-\f{|\alpha|}2+\f18}W(t)^{1-\f1{4\alpha}}\|(u_0,\omega_0)\|_{E_0^p}^{\f1{4\alpha}},
\end{align*}
and\begin{align*}
\|\omega\|_{L^\infty_T\dot B^{\f3p-\f14}_{p,\infty}}&\le
C\|\omega\|_{L^\infty_T{\dot B^{\f3p-1}_{p,\infty}}}^{\f14}\|\omega\|_{L^\infty_T{\dot B^{\f3p}_{p,\infty}}}^{\f34}\le
C\|\omega\|_{L^\infty_T{\dot B^{\f3p-1}_{p,\infty}}}^{\f14}\|D_x \omega\|_{L^\infty_T{\dot B^{\f3p-1}_{p,\infty}}}^{\f34}
\\&\le Ct^{-\f38}W(t)^{\f34}\|(u_0,\omega_0)\|_{E_0^p}^{\f1{4}}.
\end{align*}
The term of $F$ is done in the same way. Thus
\begin{align}\label{eq:decay:II}
\sup_{j\in\Z}2^{j(\f3p-1)}t^{\frac{|\alpha|}{2}}\mathcal{II}&\le CW(t)^{\f74-\f1{4\alpha}}\|(u_0,\omega_0)\|_{E_0^p}^{\f1{4}+\f1{4\alpha}}.
\end{align}

For the estimate of $\omega_d$, we localize the third equation of \eqref{eq:micropolartransform}, then taking $D^\alpha_x$ on the
localized equation yields
\begin{align*}
\p_t\Delta_jD^\alpha_x\omega_d-2\Delta\Delta_jD^\alpha_x\omega_d+2\Delta_jD^\alpha_x\omega_d
=-\Lambda^{-1}\dv D^\alpha_x\Delta_j\big(u\cdot\na\omega\big).
\end{align*}
Multiplying  by
$p|\Delta_jD^\alpha_x\omega_d|^{p-2}\Delta_jD^\alpha_x\omega_d$ and
integrating with respect to $x$ yield that
\begin{align*}
\f d{dt}&\|\Delta_jD^\alpha_x\omega_d\|_{L^p}^p+2p\int_{\R^3}(-\Delta)\Delta_jD^\alpha_x\omega_d|\Delta_jD^\alpha_x\omega_d|^{p-2}\Delta_jD^\alpha_x\omega_d{\rm d}x
\nonumber\\&+2p\int_{\R^3}|\Delta_jD^\alpha_x\omega_d|^p{\rm d}x
=-p\int_{\R^3}\Lambda^{-1}\dv D^\alpha_x\Delta_j\big(u\cdot\na\omega\big)|\Delta_jD^\alpha_x\omega_d|^{p-2}\Delta_jD^\alpha_x\omega_d{\rm d}x.
\end{align*}
Using  Lemma \ref{Lem:positiveinequality} produces that
\begin{align*}
&\f d{dt}\|\Delta_jD^\alpha_x\omega_d\|_{L^p}^p+c_p(2^{2j}+1)\|\Delta_jD^\alpha_x\omega_d\|_{L^p}^p\le C\|\Delta_jD^\alpha_x(u\cdot\na\omega)\|_{L^p}
\|\Delta_jD^\alpha_x\omega_d\|_{L^p}^{p-1}.
\end{align*}
This together with Gronwall's inequality implies that
\begin{align*}\label{}
\|\Delta_jD^\alpha_x\omega_d\|_{L^p}\le & e^{-c_p t(2^{2j}+1)}\|\Delta_jD^\alpha_x\omega_{0,d}\|_{L^p}+\mathcal{III},\end{align*}where
\begin{align*}\mathcal{III}=C\int_0^te^{-c_p (t-\tau)2^{2j}}e^{-(t-\tau)}\|D^\alpha_x\Delta_j(u\cdot\na\omega)(\tau)\|_{L^p}{\rm d}\tau.
\end{align*}
Using Lemma \ref{Lem:Bernstein} and \eqref{eq:expotentialrelation}, we obtain
\begin{align*}\label{}
\mathcal{III}\le& C\int_0^{t/2}e^{-\tilde{c}_p (t-\tau)2^{2j}}e^{-(t-\tau)}(t-\tau)^{-|\al|/2}\|\Delta_j(u\cdot\na\omega)(\tau)\|_{L^p}{\rm d}\tau\nonumber\\
&+C\int_{t/2}^te^{-\tilde{c}_p (t-\tau)2^{2j}}e^{-(t-\tau)}(t-\tau)^{-1/2}\|D^{\alpha-1}_x\Delta_j(u\cdot\na\omega)(\tau)\|_{L^p}{\rm d}\tau.
\end{align*}
The first term is treated as $\mathcal{I}$, the second term is treated as $\mathcal{II}$, then
\begin{align}\label{eq:decay:III}
\sup_{j\in\Z}2^{j(\f3p-1)}t^{\frac{|\alpha|}{2}}\mathcal{III}&\le CW(t)^{\f74-\f1{4\alpha}}\|(u_0,\omega_0)\|_{E_0^p}^{\f1{4}+\f1{4\alpha}}.
\end{align}
Combining \eqref{eq:decay} with \eqref{eq:decay:I}--\eqref{eq:decay:III}, we have
\begin{align*}
A(t)\le \|(u_0,\omega_0)\|_{E_0^p}+C\|(u_0,\omega_0)\|_{E_0^p}^2+CW(t)^{\f74-\f1{4\alpha}}\|(u_0,\omega_0)\|_{E_0^p}^{\f1{4}+\f1{4\alpha}}.
\end{align*}
Then by the continuous induction, we have $W(t)\le 2CE$. This complete the proof of Theorem \eqref{Thm:decaytheorem}.

\section{Acknowledgments}
Q. Chen and C.
Miao  were   supported by the NSF of China  grants 10701012 and 10725102 respectively).

\end{document}